\documentclass[final]{siamltex}
\usepackage{graphicx}
\usepackage[english]{babel}

\def\beq{\begin{equation}}
\def\eeq{\end{equation}}
\def\ce{{\cal E}}
\def\cf{{\cal F}}
\def\cl{{\cal L}}
\def\cm{{\cal M}}
\def\pn{{\cal P}}
\def\qn{{\cal Q}}

\def\cn{{\cal N}}
\newcommand{\prefix}[1]{{}^{#1}}

\title{Generalized Selective Modal Analysis }
\author{Juli\'an Barqu\'{\i}n  
        \thanks{Instituto de Investigaci\'on Tecnol\'ogica,
        Universidad Pontificia Comillas de Madrid,
        C/ Alberto Aguilera, 23,
        E-28015 Madrid, Spain ({\tt barquin@iit.upco.es}).}}

\begin{document}

\pagestyle{myheadings}
\thispagestyle{plain}
\markboth{J. BARQUIN}{GENERALIZED SELECTIVE MODAL ANALYSIS}

\maketitle

\begin{abstract}
A new approach which generalizes the Selective Modal Analyis
(SMA) and algorithms based upon it for solving the generalized
eigenvalue problem is described. This approach allows for the
systematic consideration of physical properties of the system
under study. Two small application cases
demonstrate the capabilities of the proposed approach.
\end{abstract}


%\begin{AMS}
%65F15, 15A18
%\end{AMS}

\section{Review of SMA}

Selective Modal Analysis (SMA) is a
physically motivated framework for 
understanding, simplifying and 
analyzing complex linear time invariant
models of dynamic systems. SMA can focus
on selected portions of the structure and
behaviour of the system \cite{Perez, Auto}.

In many physical systems, it can be
readily recognized that some set of modes
is asociated to a certain set of variables.
For instance:
in electrical power systems, the
          electromechanical oscillations are
          associated to the machines rotors's
          angle and speed.
Very often, this association is used,
explictly or implicitly, to simplify the
mathematical models of the system under
study.
SMA aims to exploit this relationship in a 
systematic and rigorous way.

So, let us assume that it is desired to
analyze the dynamic system:

\beq
\dot{x} = Ax
\label{uno}
\eeq

The SMA approach is to classify the
$x$ components in relevant ($r$) and
less-relevant ($z$) components. So,
possibly after a trivial reordering, it can
be written:

\beq
x = \left[ \begin{array}{c}
     r \\ z \end{array} \right]
\eeq

Therefore, equation (~\ref{uno}~) can be written
as:

\beq
\left[ \begin{array}{c}
    \dot{ r} \\ \dot{z} 
\end{array} \right]
= 
\left[ \begin{array}{cc}
    A_{rr} & A_{rz} \\ A_{zr} & A_{zz}
\end{array} \right]
\left[ \begin{array}{c}
     r \\ z 
\end{array} \right]
\label{tres}
\eeq

Let us assume that there is interest in
computing and eigenvalue $\lambda$, and
its left and right eigenvectors $v$ and $w$.

\begin{eqnarray}
\lambda v & = & A v \\
w^T\lambda & = & w^TA
\end{eqnarray}

The eigenvectors can be partitioned 
analogously as the states $x$:

\beq
v = \left[ \begin{array}{c}
      v_r \\ v_z \end{array}\right] 
\;\;\;
w = \left[ \begin{array}{c}
      w_r \\ w_z \end{array}\right] 
\label{tresb}
\eeq

It is easy to check that it must be fulfilled that:

\begin{eqnarray}
\lambda v_r & = & \left( A_{rr} + A_{rz} (\lambda - A_{zz})^{-1} A_{zr} \right) v_r \\
w_r^T\lambda & = & w_r^T \left( A_{rr} + A_{rz} (\lambda - A_{zz})^{-1} A_{zr} \right) 
\end{eqnarray}

So, the interesting eigenvalue $\lambda$ is in the spectrum
of the matrix $A_{rr} + A_{rz} (\lambda - A_{zz})^{-1} A_{zr}$. 
On the other hand, if this mode is strongly correlated
with the relevant variables $r$, it should be expected
that the spectrum of $A_{rr}$ contains an eigenvalue
quite similar to $\lambda$ and, therefore, that the matrix
$A_{rz} (\lambda - A_{zz})^{-1} A_{zr}$ perturbs
sligthly the desired mode. That suggests the
following algorithm:

\vspace{3mm}
\begin{center}
\fbox{
\begin{minipage}{100cm}
\begin{tabbing}
{\bf Algorithm 1} \\
\\
{\bf Input:} $A_{rr},A_{rz},A_{zr},A_{zz}$. \\
\\
{\bf Output:} $\lambda, v_r, w_r$. \\
\\
{\bf 1.} \= Perform the eigenanalysis of $A_{rr}$, \\
{\bf 2.} \> Select the interesting mode ${}^0\lambda, {}^0v_r, {}^0w_r$, \\
{\bf 3.} \> {\bf for} $j=1,2,3, \ldots$ until convergence, \\
\> {\bf 3.1.} \= Compute $H(\prefix{j-1}\lambda) = A_{rz} (\prefix{j-1}\lambda - A_{zz})^{-1} A_{zr}$. \\
\> {\bf 3.2.} \> Perform the eigenanalysis of $A_{rr} + H(\prefix{j-1}\lambda)$, \\
\> {\bf 3.3.} \> Select the interesting mode $\prefix{j}\lambda, \prefix{j}v_r, \prefix{j}w_r$, \\
{\bf 4.} \> {\bf end}
\end{tabbing}
\end{minipage}
}
\end{center}
\vspace{3mm}

 The convergence properties of this algorithm have been
studied in \cite{TesisPerez}. The convergence is controlled
by the eigenvectors. Especifically,
the algorithm locally converges if and only if it is fulffiled:

\beq
\mid  \rho\mid  = \left| \frac{w_r^Tv_r}{w_z^Tv_z} \right| > 1 \label{partc}
\eeq

The number $\rho$ is called the participation ratio.
Notice that, for any eigenvalue, this ratio only depends
in the way that the variables have been partitioned
in relevant and less relevant.

 The former algorithm can be generalized to search for
several eigenvalues. The interesting eigenvalues shall
be collected in a diagonal matrix $\Lambda$, and the
eigenvectors in matrices $V$ and $W$:

\beq
\Lambda = \left[ \begin{array}{cccc}
                               \lambda_1 & 0 & \ldots & 0 \\
                                0  & \lambda_2 & \ldots &0 \\
                               \vdots & \vdots & \ddots & \vdots \\
                              0 & 0 & \ldots & \lambda_n
                    \end{array} \right]
\;\;
V = \left[ v_1 v_2 \ldots v_n \right]
\;\;
W = \left[ w_1 w_2 \ldots w_n \right]
\eeq

The matrices $V$ and $W$ can be also partitioned
in relevant and less relevant parts:

\beq
V = \left[ \begin{array}{c} V_r \\ V_z \end{array} \right]
\;\;
W = \left[ \begin{array}{c} W_r \\ W_z \end{array} \right]
\eeq

So, the following algorithm can be proposed:

\vspace{3mm}
\fbox{
\begin{minipage}{100cm}
\begin{tabbing}
{\bf Algorithm 2} \\
\\
{\bf Input:} $A_{rr},A_{rz},A_{zr},A_{zz}$. \\
\\
{\bf Output:} $\Lambda, V_r, W_r$. \\
\\
{\bf 1.} \= Perform the eigenanalysis of $A_{rr}$, \\
{\bf 2.} \> Select the interesting modes ${}^0\Lambda, {}^0V_r, {}^0W_r$, \\
{\bf 3.} \> {\bf for} $j=1,2,3, \ldots$ until convergence, \\
\> {\bf 3.1.} \= Compute  a matrix $\prefix{j}M$
   which fulfills \\
\> \>
$\prefix{j}M\prefix{j-1}V_r = 
\left[
H(\prefix{j-1}\lambda_1) \prefix{j-1}v_{1r},
H(\prefix{j-1}\lambda_2) \prefix{j-1}v_{2r},
\ldots
H(\prefix{j-1}\lambda_n) \prefix{j-1}v_{nr}
\right]$ \\
\> \> where $H(\lambda) =
A_{rz} (\lambda - A_{zz})^{-1} A_{zr}$, \\
\> {\bf 3.2.} \> Perform the eigenanalysis of $A_{rr} + \prefix{j}M$, \\
\> {\bf 3.3.} \> Select the interesting modes $\prefix{j}\Lambda, \prefix{j}V_r, \prefix{j}W_r$, \\
{\bf 4.} \> {\bf end}
\end{tabbing}
\end{minipage}
}
\vspace{3mm}

The convergence conditions of algorithm 2 are also 
studied in \cite{TesisPerez}, although they are considerably
more involved than those of algorithm 1. However,
the computational experience shows that the convergence
is good if the participation ratios of the interesting modes
are high.

From the point of view of the computational effort,
the most demanding task of both algorithms is the
step 3.1, which requires to solve linear system
involving the matrix $\lambda - A_{zz}$. Most 
SMA applications have been developped for
the study of electric power systems \cite{SMAx},
where special techniques based upon the peculiar
characteristics of these systemas have been 
used to perform efficiently this task.

In addition to algorithms 1 and 2, there are a
number of related ones which considerably
improve their speed and robustness \cite{Sancha,Rouco}.
However, these algorithms are outside of the scope of
this paper.

\section{Generalized SMA}

Although SMA has been succesfully used in
applicatiosn, it has some drawbacks:

\begin{remunerate}

\item Very often, the computation of the
           desired modes begin with the study of
           a simplified model, in order to gain an
           insight on the interesting eigenstructure.
           Although SMA allows to use subsequently
           the information on the relevant variables,
           it does not make use of all the information
           obtained with the simplified model.

\item There are some problems, specially
           in continous media, where it is known
          the overall shape of the desired modes,
          but it is impossible to adscribe them to
         a small number of system
          variables.

\end{remunerate}

Both kind of problems are treated in this paper
examples. The aim of this section is to generalize
the classical SMA theory in order to deal with these
problems.

So, the problem to solve is the eigenvalue problem:

\begin{eqnarray}
\lambda E v & = & A v \label{eig1} \\
w^{\dagger}E\lambda & = & w^{\dagger}A
\end{eqnarray}

$A$ and $E$ are $m \times m$ real matrices.
The matrix $E$ is a symmetric, possibly singular,
projection matrix:

\beq
E^2 = E = E^T = E^{\dagger}
\eeq

The superscript $T$ denotes the transpose and
$\dagger$ the hermitian conjugate.
It is assumed that the 
right eigenvector $v$ approximately
lies in the subspace spanned by 
$\{ e_1, e_2, \ldots, e_n \}$.
Usually, $n \ll m$.
Then, it is defined the matrix

\beq
\ce = \left[ e_1, e_2, \ldots, e_n \right]
\eeq

 In similar way, the left eigenvector $w$ is
assumed to yield, approximately, in the subspace
spanned by
$\{ f_1,f_2, \ldots, f_n \}$.
So, it is defined the matrix

\beq
\cf = \left[ f_1,f_2, \ldots, f_n \right]
\eeq

 Besides, the $e_i$ and $f_j$ basis are normalized
in order to fulfill the equation:

\beq
\cf ^{\dagger}E\ce = I_n
\label{norm1}
\eeq

$I_n$ is the $n \times n$ identity matrix.
This equation can be enforced so long as no vector
generated by the basis $e_i$ or the basis
$f_j$
is included in the kernel of $E$. This condition
shall be assumed in the sequel.

Then, the eigenvectors $v$ nd $w$ can be written as:

\begin{eqnarray}
 v & = & \ce \alpha + z \label{deco1} \\
\cf ^{\dagger} E z & = & 0 \label{deco2} \\
 w & = & \cf \beta + y \label{deco3} \\
\ce ^{\dagger} E y & = & 0 \label{deco4}
\end{eqnarray}
 
 It is easy to show that the above decomposition exists and
is unique.  Then, after the algebraic manipulations
shown in the appendix A, it is found that:

\beq
\lambda \alpha =
A_{rr}\alpha + H(\lambda) \alpha
\eeq

\noindent where

\begin{eqnarray}
A_{rr} & = & \cf^{\dagger} A \ce \label{T1}  \\
H(\lambda) & = &
\cf^{\dagger} A
\pn\left\{ \lambda E -
A + \left[ A, \qn \right]_+
\right\}^{-1}\pn A\ce 
\label{T2}
\end{eqnarray}

The matrices $\qn$ and $\pn$ are idempotent matrices
defined by

\beq
\pn = I_m - E\ce \cf^{\dagger}E = I_m - \qn
\label{defp}
\eeq

\noindent and $\left[ \qn, A \right]_+$ is the anti-commutator:

\beq
\left[ \qn, A \right]_+ = \qn A + A \qn
\eeq

These formulae are the basic ones in selective modal analysis,
and can be considered as a generalization of ``classical'' SMA,
as shown in the appendix B.
It is also easy to check that $\beta$ is the left eigenvector
of $A_{rr} + H(\lambda)$.
Notice that the dimension of the matrices $A_{rr}$ and 
$H(\lambda)$ is $n \ll m$.

It is also noteworhy that

\begin{eqnarray}
z & = & 
\pn\left\{ \lambda E -
A + \left[ A, \qn \right]_+
\right\}^{-1}\pn A\ce \alpha  \\
y^{\dagger} & = &
\beta^{\dagger} \cf^{\dagger}  A \pn
\left\{ \lambda E -
A + \left[ A, \qn \right]_+ 
\right\}^{-1} \pn 
\end{eqnarray}

\section{Algorithms based on  $H(\lambda)$ computation}

The aim of this section is to apply the results of the
former section in order to obtain workable algorithms,
and their convergence conditions. In order to simplify
the notation, let us denote by $N(\lambda,\qn)$ 
the matrix

\beq
N(\lambda,\qn)  = \left\{ \lambda E - A + [A,\qn]_+\right\}^{-1} 
\eeq

\noindent so that

\beq
H(\lambda) = 
\cf^{\dagger} A \pn 
\left\{ \prefix{j-1}\lambda E -
A + \left[ A, \qn \right]_+
\right\}^{-1} \pn A \ce = 
\cf^{\dagger} A \pn 
N(\lambda,\qn)
\pn A \ce
\eeq

\subsection{Linear algorithm}

Specifically, let us consider the following generalization
of algorithm 1:

\vspace{3mm}
\begin{center}
\fbox{
\begin{minipage}{100cm}
\begin{tabbing}
{\bf Algorithm 3} \\
\\
{\bf Input:} $E,A, \ce, \cf$. \\
\\
{\bf Output:} $\lambda , \alpha , \beta$. \\
\\
{\bf 1.} \= Form $A_{rr} = \cf^{\dagger}A\ce$, and perform the eigenanalysis of $A_{rr}$, \\
{\bf 2.} \> Select the interesting mode ${}^0\lambda, {}^0\alpha, {}^0\beta$, \\
{\bf 3.} \> {\bf for} $j=1,2,3, \ldots$ until convergence, \\
\> {\bf 3.1.} \= Compute $H(\prefix{j-1}\lambda) = 
\cf^{\dagger} A
\pn N(\prefix{j-1}\lambda,\qn )   \pn A\ce$. \\
\> {\bf 3.2.} \> Perform the eigenanalysis of $A_{rr} + H(\prefix{j-1}\lambda)$, \\
\> {\bf 3.3.} \> Select the interesting mode $\prefix{j}\lambda, \prefix{j}\alpha, \prefix{j}\beta$, \\
{\bf 4.} \> {\bf end}
\end{tabbing}
\end{minipage}
}
\end{center}
\vspace{3mm}

The following theorem states the conditions for the local
convergence of Algorithm 3:

\begin{theorem}
\label{teo1}
Given an eigenvalue $\lambda$ of the pair $(E,A)$
with associated right and left eigenvectors $v$ and $w$,
there is $\delta > 0$ such that if 
$\| \prefix{0}\lambda - \lambda \| < \delta$, algorithm
3 converges to the eigenvalue $\lambda$ whenever 

\beq
| \rho | = 
\left| \frac{w^{\dagger} \qn v}{w^{\dagger} (E-\qn) v} \right| > 1
\eeq

Furtthermore, the error $\prefix{j}\epsilon =
\prefix{j}\lambda - \lambda$ fulfills:

\beq
 \prefix{j}\epsilon  = \rho \prefix{j-1}\epsilon + o(\prefix{j-1}\epsilon)
\eeq

\end{theorem}

\begin{proof}
The proof is given in the appendix C.
\end{proof}

It is easy to check that $\rho$ just defined
is, in the ``classical'' case, the same $\rho$
defined in equation (\ref{partc}).

\subsection{Superlinear algorithm}

If, in algorithm 3,
the right eigenvector $v$
(respectively the left eigenvector $w$)
is contained in the span of $\ce$
($\cf$), then
$z = 0$ ($y = 0$) and $\rho \rightarrow \infty$.
So, the algorithm could be speeded up if
the matrices $\ce$ and $\cf$
are pdated in order that their
span contains the last approximation
to $v$ and $w$: $\prefix{j}v$ and $\prefix{j}w$.
Therefore, the following algorithm is
proposed:

\vspace{3mm}
\begin{center}
\fbox{
\begin{minipage}{100cm}
\begin{tabbing}
{\bf Algorithm 4} \\
\\
{\bf Input:} $E,A, {}^0\ce, {}^0\cf$. \\
\\
{\bf Output:} $\lambda , \alpha , \beta, v, w$. \\
\\
{\bf 1.} \= Form $\prefix{0}A_{rr} = {}^0\cf^{\dagger}A\;{}^0\ce$, 
           and perform the eigenanalysis of $A_{rr}$, \\
{\bf 2.} \> Select the interesting mode ${}^0\lambda, {}^0\alpha, {}^0\beta$, \\
{\bf 3.} \> {\bf for} $j=1,2,3, \ldots$ until convergence, \\
\> {\bf 3.1.} \= Compute $H(\prefix{j-1}\lambda) = $ \\
\> \>
$\prefix{j-1}\cf^{\dagger} A\;
\prefix{j-1}\pn
N(\prefix{j-1}\lambda,\prefix{j-1}\qn)
\;\prefix{j-1}\pn A\;\prefix{j-1}\ce$. \\
\> {\bf 3.2.} \> Perform the eigenanalysis of $\prefix{j-1}A_{rr} + H(\prefix{j-1}\lambda)$, \\
\> {\bf 3.3.} \> Select the interesting mode $\prefix{j}\lambda, 
                   \prefix{j}\alpha, \prefix{j}\beta$, \\
\> {\bf 3.4.} \> Compute 
                  $\prefix{j}z = \prefix{j-1}\pn N(\prefix{j-1}\lambda,\prefix{j-1}\qn)
                                      \;\prefix{j-1}\pn A \prefix{j-1}\ce \prefix{j}\alpha$ \\
\> \>                   and 
                   $\prefix{j}y^{\dagger} = \prefix{j}\beta^{\dagger}\prefix{j-1}\cf^{\dagger}
                                     A\prefix{j-1}\pn N(\prefix{j-1}\lambda,\prefix{j-1}\qn)
                                      \;\prefix{j-1}\pn$, \\
\> {\bf 3.5.} \> Compute
                  $\prefix{j}v = \prefix{j-1}\ce \prefix{j}\alpha + \prefix{j}z$ and
                   $\prefix{j}w^{\dagger} = \prefix{j}\beta^{\dagger} \prefix{j-1}\cf^{\dagger}
                                         + \prefix{j}y^{\dagger}$, \\                                                                                                                                 
\> {\bf 3.6.} \> Update $\prefix{j}\ce$, $\prefix{j}\cf$ in such a way that \\
\> \> $\prefix{j}v \in {\rm span}(\prefix{j}\ce), \prefix{j}w \in {\rm span}(\prefix{j}\cf)$, \\
\> {\bf 3.7} \> Form $\prefix{j}A_{rr} = {}^j\cf^{\dagger}A\;{}^j\ce$, \\
{\bf 4.} \> {\bf end}
\end{tabbing}
\end{minipage}
}
\end{center}
\vspace{3mm}

A particular case of algorithm 4 is when
the matrices $\prefix{j}\ce$ and
$\prefix{j}\cf$ are vectors. Then,
these matrices are esentially the
estimated eigenvectors. The local 
convergence properties of the algorithm
are, in this case, given by the following theorem:

\begin{theorem}
\label{teo2}
Given an eigenvalue $\lambda$ of the pair
$(E,A)$ with associated right and left
eigenvectors $v$ and $w$, 
if $\| N(\lambda, Evw^{\dagger}E)\| < \infty$,
there is a neighborhood of $\lambda, v, w$,
such that if $\prefix{0}\lambda, \prefix{0}v,
\prefix{0}w$ belong to it, algorithm 4
converges. Furthermore, it is fullfilled
that asymptotically there is a constant
$K$ such that

\beq
| \prefix{j}\epsilon | \leq  K | \prefix{j-1}\epsilon |^{1 + \sqrt{2}}
\eeq

\end{theorem}

\begin{proof}
The proof is given in the appendix \ref{aped}.
\end{proof}

In the general case, whenever $\prefix{j}\ce$
and $\prefix{j}\cf$ are not vectors, it is expected
that algorithm 4 converges at least so fast.
This is because the relevant subspace is bigger,
so that the approximation to the eigenstructure
can not be worst.

\subsection{$\ce$ and $\cf$ selection}

The computation of $H(\lambda)$ requieres
to make a selection of the matrices $\ce$ and
$\cf$. The following theorem can be used
for this task:

\begin{theorem}
\label{teoce}

The matrix $A_{rr} + H(\lambda)$ is invariant
under the transformations 
$\ce \leftarrow \ce + (I_m - E) \cl$ or
$\cf \leftarrow \cf + (I_m - E) \cm$,
where $\cl$ or $\cm$ are arbitrary matrices
of the same dimension than $\ce$ or $\cf$.

\end{theorem}

\begin{proof}

The proof is provided in appendix
\ref{apecc}

\end{proof}

 In many cases the $E$ matrix can be written as:

\beq
E = \left[ \begin{array}{cc}
       I_r & 0 \\ 0 & 0
      \end{array} \right]
\eeq

$I_r$ is a $r$-dimensional identity matrix.
Therefore, the vectors and matrices can
be partitioned in dynamic and static parts.
Specifically,

\begin{eqnarray}
\ce & = & \left[ \begin{array}{c}
                       \ce_d \\ \ce_s
               \end{array} \right] \\
\cf & = & \left[ \begin{array}{c}
                       \cf_d \\ \cf_s
               \end{array} \right] 
\end{eqnarray}

The invariance of the matrix $A_{rr} + H(\lambda )$
under the considered transformations
means that the value of the static components
$\ce_s$ and $\cf_s$ is irrelevant in order
to compute this matrix. 

So, referring to algorithm 4, there are 
at least two possibilities:

\begin{remunerate}

\item To keep the whole eigenvector in step 3.4:
update $\prefix{j}\ce$, $\prefix{j}\cf$ in such a way that 
$\prefix{j}v \in {\rm span}(\prefix{j}\ce), \prefix{j}w \in {\rm span}(\prefix{j}\cf)$.
Then, as the algorithm converges, the matrix $H(\lambda) \rightarrow 0$.

\item To update the matrix as above, but 
the static components, which are zeroed ($\ce_s=0$ and $\cf_s=0$).
As the matrix converges, the matrix $H(\lambda)$ converges to a non-zero
value.

\end{remunerate}

The second possibility can be useful in order to minimize the
numerical effort. 

\section{$H(\lambda )$ computation}

From the computational point of view,
the most demanding task of the algorithm is
the computation of the matrix $H(\lambda)$ or
of the matrix $A_{rr} + H(\lambda )$. The purporse
of this section is to propose algorithms to deal
efficiently with this task.

\subsection{$H(\lambda)$ computation using the Shermann-Morrison lemma}

$H(\lambda )$ can be written as:

\beq
H(\lambda) = \cf^{\dagger} A \pn \left\{
              \lambda E - A + 
              \left[A,\qn\right]_+
              \right\} ^{-1} \pn A\ce
\eeq

The basic problem is related to the matrix

\beq
N(\lambda) = 
\left\{
              \lambda E -  A + 
              \left[A,\qn\right]_+
              \right\} ^{-1}
\eeq

A problem is that, generally, the matrix
$\left[A,\qn \right]_+$ is not sparse. 
However, it is possible to obtain an
expression of $N(\lambda)$, which allows computations
by using only sparse matrices, by means of the
Shermann-Morrison lemma. So,

\begin{eqnarray}
N(\lambda) & = & 
\left\{
              \lambda E -  A + 
              \left[A,\qn \right]_+ 
              \right\} ^{-1} \\
& = &
\left\{
              \lambda E - A + 
              \left[A,E\ce\cf^{\dagger}E \right]_+
              \right\} ^{-1} \\
& = &
\left\{
              \lambda E - A + 
              \left( AE\ce \right)\left(\cf^{\dagger}E\right)
            + \left(E\ce\right) \left(\cf^{\dagger}EA\right)
              \right\} ^{-1} \\
& = &
\left\{
              \lambda E - A +
              \left( AE\ce \right)\left(\cf^{\dagger}E\right)
            + \left(E\ce\right) \left(\cf^{\dagger}EA\right)
              \right\} ^{-1} \\
& = &
\left\{
              \lambda E - A + \eta\phi^{\dagger} - \eta\phi^{\dagger} + 
\right.
\nonumber \\ & & \left.
              \left( AE\ce \right)\left(\cf^{\dagger}E\right)
            + \left(E\ce\right) \left(\cf^{\dagger}EA\right)
              \right\} ^{-1} \\
& = &
\left\{
              \lambda E - A + \eta\phi^{\dagger} - 
              \left[ \begin{array}{ccc}
              \eta & -AE\ce & -E\ce \end{array}
              \right]
              \left[ \begin{array}{c}
              \phi^{\dagger} \\ \cf^{\dagger} E \\ \cf^{\dagger} EA 
              \end{array} \right] 
              \right\} ^{-1} 
\end{eqnarray}

$\eta$ and $\phi$ are two sparse vectors which make sure
that the matrix $\lambda E - A + \eta\phi^{\dagger}$ is
regular even if $\lambda$ is an eigenvalue of the pair
$(A,E)$. This is going to happen when the SMA algorithm
converges. Now, from the Shermann-Morrison lemma:

\begin{eqnarray}
N(\lambda) & = & 
\left\{ \lambda E - A + \eta\phi^{\dagger} \right\} ^{-1} +  \nonumber \\
& &
\left\{ \lambda E - A + \eta\phi^{\dagger} \right\} ^{-1} 
              \left[ \begin{array}{ccc}
              \eta & -AE\ce & -E\ce \end{array}
              \right] \nonumber \\ & &
\left(
I_{2n+1} -
              \left[ \begin{array}{c}
              \phi^{\dagger} \\ \cf^{\dagger} E \\ \cf^{\dagger} EA 
              \end{array} \right] 
\left\{ \lambda E - A + \eta\phi^{\dagger} \right\} ^{-1} 
              \left[ \begin{array}{ccc}
              \eta & -AE\ce & -E\ce \end{array}
              \right]
\right)^{-1} \nonumber \\ & &
              \left[ \begin{array}{c}
              \phi^{\dagger} \\ \cf^{\dagger} E \\ \cf^{\dagger} EA 
              \end{array} \right] 
\left\{ \lambda E - A + \eta\phi^{\dagger} \right\} ^{-1} 
\label{sh-mo}
\end{eqnarray}

The number $2n+1$ is usually small. So, it is only required
to know the LU factorization of a filled matrix of small
($2n+1$) dimension and of the sparse matrix
$\lambda E - A + \eta\phi^{\dagger} $.
 
\subsection{Composite models}

In power systems analysis, the system to analyze
is a set of dynamical subsystems connected through a
static relationship. Specifically, there are $l$ subsytems

\begin{eqnarray}
E_k \dot{x}_{Mk} & = & A_k x_{Mk} + B_k x_{Ik}  \label{compo1} \\
x_{Ok} & = & C_k x_{Mk} + D_k x_{Ik} \;\;\;\;\; k = 1,\ldots,l \label{compo2}
\end{eqnarray}

The variables $x_{Mk}$ are the state variables of 
the $k$-th subsytem, $x_{Ik}$ are the input variables
and $x_{Ok}$ the output variables. It is assumed that
the number of output and input variables of each system is equal.
It is also assumed that
the matrices $E_k$ are symetric projection real matrices:

\beq
E_k = E_k^T = E_k^{\dagger} = E_k^2
\eeq

Let us define the vectors

\beq
x_M = \left[ \begin{array}{c}
            x_{M1} \\ x_{M2} \\ \vdots \\ x_{Ml}
           \end{array} \right]
\;\;
x_I = \left[ \begin{array}{c}
            x_{I1} \\ x_{I2} \\ \vdots \\ x_{Il}
           \end{array} \right]
\;\;
x_O = \left[ \begin{array}{c}
            x_{O1} \\ x_{O2} \\ \vdots \\ x_{Ol}
           \end{array} \right]
\eeq

In addition to these equations,
there is also a static interconnection:

\beq
\left[ \begin{array}{cc}
J_{11} & J_{12} \\
J_{21} & J_{22}
\end{array} \right]
\left[ \begin{array}{c}
x_I \\ x_A
\end{array} \right] 
=
\left[ \begin{array}{c}
x_O \\ 0
\end{array} \right] 
\label{compo3}
\eeq

$x_A$ is a set of additional algebraic variables.

For this kind of systems, it is convenient
to consider the following $\ce$ and $\cf$
matrices:

\beq
\ce = \left[ \begin{array}{c}
                    \ce_M \\ 0 \\ 0 \\ 0
          \end{array} \right]
\;\;\;\;
\cf = \left[ \begin{array}{c}
                    \cf_M \\ 0 \\ 0 \\ 0
          \end{array} \right]
\eeq

\noindent which mimics the $x$ structure. Furthermore,
$\ce_M$ and $\cf_M$ are defined on a subsytem basis:

\begin{eqnarray}
\ce_M & = & \left[ \begin{array}{cccc}
               \ce_{M1} & 0 & \ldots & 0 \\
               0 & \ce_{M2} & \ldots & 0 \\
               \vdots & \vdots & \ddots & \vdots \\
              0 & 0 & \ldots & \ce_{Ml} 
             \end{array} \right] = {\rm diag} (\ce_{M1}, \ce_{M2}, \ldots, \ce_{Ml} ) \\
 \cf_M & = & \left[ \begin{array}{cccc}
               \cf_{M1} & 0 & \ldots & 0 \\
               0 & \cf_{M2} & \ldots & 0 \\
               \vdots & \vdots & \ddots & \vdots \\
              0 & 0 & \ldots & \cf_{Ml} 
             \end{array} \right] = {\rm diag} (\cf_{M1}, \cf_{M2}, \ldots, \cf_{Ml} ) 
\end{eqnarray}

Then, as proved in the appendix \ref{apef},
the matrix $H(\lambda)$ can be computed as

\beq
H(\lambda) =
H_A(\lambda) +
\left[
\left( B_r + H_B(\lambda) \right) \; 0 \right]
\left[ \begin{array}{cc}
J_{11}-(D+H_D(\lambda)) & J_{12} \\
J_{21} & J_{22}
\end{array} \right]^{-1}
\left[ \begin{array}{c} 
C_r + H_C(\lambda) \\ 0 
\end{array} \right]
\eeq

\noindent where all the matrices are computed in a subsytem basis:

\begin{eqnarray}
A_r & = & {\rm diag} (A_{r1} \ldots A_{rl}) \\
B_r & = & {\rm diag} (B_{r1} \ldots B_{rl}) \\
C_r & = & {\rm diag} (C_{r1} \ldots C_{rl}) \\
D & = & {\rm diag} (D_{1} \ldots D_{l}) \\
H_{A} & = & {\rm diag} (H_{A1} \ldots H_{Al}) \\
H_{B} & = & {\rm diag} (H_{B1} \ldots H_{Bl}) \\
H_{C} & = & {\rm diag} (H_{C1} \ldots H_{Cl}) \\
H_{D} & = & {\rm diag} (H_{D1} \ldots H_{Dl}) 
\end{eqnarray}

\section{Direct algorithms}

As said above, the most difficult task in order to
apply SMA is the computation of the matrix
$H(\lambda)$. However, the proposed algorithms
can be formulated without needing to compute
this matrix. The purporse of this section is to explain
the way of doing it.

\subsection{Single eigenvalue algorithms}

 The basic SMA formulae are:

\begin{eqnarray}
\prefix{j} \lambda \prefix{j} \alpha & = &
\left( A_{rr} + H(\prefix{j-1}\lambda) \right) 
 \prefix{j} \alpha \label{basica1} \\
\prefix{j} \beta^{\dagger} \prefix{j} \lambda & = &
\prefix{j} \beta^{\dagger} 
\left( A_{rr} + H(\prefix{j-1}\lambda) \right) 
\label{basica2} \\
\prefix{j}z & = & 
\pn\left\{ \prefix{j-1}\lambda E -
A + \left[ A, \qn \right]_+
\right\}^{-1}\pn A\ce \prefix{j}\alpha  \label{basica3} \\
\prefix{j}y^{\dagger} & = &
\beta^{\dagger} \cf^{\dagger}  A \pn
\left\{ \prefix{j-1}\lambda E -
A + \left[ A, \qn \right]_+ 
\right\}^{-1} \pn \label{basica4}
\end{eqnarray}

\noindent as shown in prvious sections. Let us also, as above,
define the vectors

\begin{eqnarray}
\prefix{j}v & = & \ce \prefix{j} \alpha + \prefix{j}z \label{basica5} \\
\prefix{j}w & = & \cf \prefix{j} \beta + \prefix{j} y \label{basica6}
\end{eqnarray}

Then, the following theorem can be stated:

\begin{theorem}
\label{teof}
Let us assume that equations (\ref{basica1}-\ref{basica4}) are fullfilled.
Then, the vectors $\prefix{j}v$ and $\prefix{j}w$ fullfill the equation

\begin{eqnarray}
\left[ A - \prefix{j-1}\lambda \left( E - \qn \right) \right] \prefix{j}v 
& = & \prefix{j}\lambda \qn \prefix{j}v \label{dedu1} \\
\prefix{j}w^{\dagger}\left[ A - \prefix{j-1}\lambda \left( E - \qn \right) \right] 
& = & \prefix{j}w^{\dagger}\qn \prefix{j}\lambda  \label{dedu2} 
\end{eqnarray}

\end{theorem}

\begin{proof}
The proof is provided in Appendix \ref{apeg}
\end{proof}

Let us write equation (\ref{dedu1}) in the following equivalent form

\beq
\left( A - \prefix{j-1}\lambda E \right)  \prefix{j}v  = 
\left(  \prefix{j}\lambda  -  \prefix{j-1}\lambda \right) \qn  \prefix{j}v \label{dedu3} \\
\eeq

As $\cf^{\dagger} E \prefix{j}z = 0$, it is obtained

\beq
\qn \prefix{j}v = E \ce \prefix{j}\alpha \label{dedu3b}
\eeq

Therefore

\beq
\left( A - \prefix{j-1}\lambda E \right)  \prefix{j}v  = 
\left(  \prefix{j}\lambda  -  \prefix{j-1}\lambda \right) E \ce  \prefix{j}\alpha  
\eeq

Let us define

\beq
\prefix{j}V = \left( A - \prefix{j-1}\lambda E \right)^{-1} E \ce \label{Mar1}
\eeq

So

\beq
\prefix{j}v = \prefix{j}V  \prefix{j}\alpha  
 \left(  \prefix{j}\lambda  -  \prefix{j-1}\lambda \right) 
\eeq

Analogously,

\begin{eqnarray}
\prefix{j}W^{\dagger} & = & \cf^{\dagger} E \left( A - \prefix{j-1}\lambda E \right)^{-1}
\label{Mar2}  \\
\prefix{j}w^{\dagger} & = &  \left(  \prefix{j}\lambda  -  \prefix{j-1}\lambda \right) 
\prefix{j}\beta^{\dagger} \prefix{j}W^{\dagger}
\end{eqnarray}

Now, let us consider the matrix
$\prefix{j}W^{\dagger} \left( A - \prefix{j-1}\lambda E \right) \prefix{j}V$.

\begin{eqnarray}
\prefix{j}W^{\dagger} \left( A - \prefix{j-1}\lambda E \right) \prefix{j}V \prefix{j}\alpha
& = &
\cf^{\dagger} E \prefix{j}v \left(  \prefix{j}\lambda  -  \prefix{j-1}\lambda \right)^{-1} \nonumber \\
& = &
\left(  \prefix{j}\lambda  -  \prefix{j-1}\lambda \right)^{-1} \prefix{j}\alpha \label{Mar3}
\end{eqnarray}

Therefore, $\prefix{j}\alpha$ is a right eigenvector of this matrix with
associated eigenvalue $\left(  \prefix{j}\lambda  -  \prefix{j-1}\lambda \right)^{-1}$.
It can be easily shown that $\prefix{j}\beta$ is the left eigenvector. Note that,
because of (\ref{Mar2}), this matrix can be also writen as $\cf^{\dagger}E \prefix{j}V$.

Equations (\ref{Mar1}), (\ref{Mar2}) y (\ref{Mar3}) provide an alternative way
of performing the SMA iteration. Specifically, algorithm 3 can be also writen as:

\vspace{3mm}
\begin{center}
\fbox{
\begin{minipage}{100cm}
\begin{tabbing}
{\bf Algorithm 5} \\
\\
{\bf Input:} $E,A, \ce, \cf$. \\
\\
{\bf Output:} $\lambda , \alpha , \beta$. \\
\\
{\bf 1.} \= Form $A_{rr} = \cf^{\dagger}A\ce$, and perform the eigenanalysis of $A_{rr}$, \\
{\bf 2.} \> Select the interesting mode ${}^0\lambda, {}^0\alpha, {}^0\beta$, \\
{\bf 3.} \> {\bf for} $j=1,2,3, \ldots$ until convergence, \\
\> {\bf 3.1.} \= Compute $\prefix{j}V = \left( A - \prefix{j-1}\lambda  E \right)^{-1}E \ce$, \\
\> {\bf 3.2.} \> Compute $\prefix{j}W^{\dagger} = \cf^{\dagger}E \left( A - \prefix{j-1}\lambda 
                         E \right)^{-1}$, \\
\> {\bf 3.3.} \> Perform the eigenanalysis of $\prefix{j}{\cal M} = \cf^{\dagger}E \prefix{j}V$, \\
\> {\bf 3.4.} \> Select the interesting mode $\prefix{j}\alpha, \prefix{j}\beta, 
                        \left( \prefix{j}\lambda - \prefix{j-1}\lambda \right)^{-1}$, \\
\> {\bf 3.5.} \> Update $\prefix{j}\lambda$, \\
{\bf 4.} \> {\bf end}
\end{tabbing}
\end{minipage}
}
\end{center}
\vspace{3mm}

As this algorithm is esentially the same one that Algorithm 3, its 
convergence conditions are the same ones. 
Note also that step 3.2 could be omitted.
Analogously,
Algorithm 4 can be written as:

\vspace{3mm}
\begin{center}
\fbox{
\begin{minipage}{100cm}
\begin{tabbing}
{\bf Algorithm 6} \\
\\
{\bf Input:} $E,A, {}^0\ce, {}^0\cf$. \\
\\
{\bf Output:} $\lambda , \alpha , \beta, v, w$. \\
\\
{\bf 1.} \= Form $\prefix{0}A_{rr} = {}^0\cf^{\dagger}A\;{}^0\ce$, 
           and perform the eigenanalysis of $\prefix{0}A_{rr}$, \\
{\bf 2.} \> Select the interesting mode ${}^0\lambda, {}^0\alpha, {}^0\beta$, \\
{\bf 3.} \> {\bf for} $j=1,2,3, \ldots$ until convergence, \\
\> {\bf 3.1.} \= Compute $\prefix{j}V = \left( A - \prefix{j-1}\lambda E \right)^{-1}E \prefix{j-1}\ce$, \\
\> {\bf 3.2.} \> Compute $\prefix{j}W^{\dagger} = 
                         \prefix{j-1}\cf^{\dagger}E \left( A - \prefix{j-1}\lambda E \right)^{-1}$, \\
\> {\bf 3.3.} \> Perform the eigenanalysis of $\prefix{j}{\cal M}=\prefix{j-1}\cf^{\dagger}E \prefix{j}V$, \\
\> {\bf 3.4.} \> Select the interesting mode $\prefix{j}\alpha, \prefix{j}\beta, 
                        \left( \prefix{j}\lambda - \prefix{j-1}\lambda \right)^{-1}$, \\
\> {\bf 3.5.} \> Update $\prefix{j}\lambda$, \\
\> {\bf 3.6.} \> Update $\prefix{j}\ce$, $\prefix{j}\cf$ in such a way that \\
\> \> $\prefix{j}V \in {\rm span}(\prefix{j}\ce), \prefix{j}W \in {\rm span}(\prefix{j}\cf)$, \\
{\bf 4.} \> {\bf end}
\end{tabbing}
\end{minipage}
}
\end{center}
\vspace{3mm}

Step 3.6 warrants that 
$\prefix{j}v \in {\rm span}(\prefix{j}\ce), \prefix{j}w \in {\rm span}(\prefix{j}\cf)$.
The simplest way to achieve that is to set $\prefix{j}\ce = \prefix{j}V$ and
$\prefix{j}\cf = \prefix{j}W$ but, possibly, a normalization constant.
As previsously as Algorithm 5, convergence results related to
Algorithm 4 can be directly applied to Algorithm 6.

\subsection{Several eigenvalues algorithms}

In this section algorithm 2 is written in the direct formulation,
and a new one which includes relevant subspaces updated is
also proposed. However, no convergence results are provided.

Let us firstly define the natrix

\beq
\prefix{j}\overline{V} = \left[ A - \prefix{j-1}\lambda \left(
                                                     E - \qn \right) \right]^{-1} E \ce \label{defV}
\eeq

It is clear, from (\ref{dedu1}) and (\ref{dedu3b}):

\beq
\prefix{j}v = \prefix{j}\lambda \prefix{j}\overline{V} \prefix{j}\alpha \label{deduA}
\eeq

Analogously, it is possible to write:

\begin{eqnarray}
\prefix{j}\overline{W}^{\dagger} & = & \cf^{\dagger}E \left[ A - \prefix{j-1}\lambda \left(
                                                     E - \qn \right) \right]^{-1} \\
\prefix{j}w^{\dagger} & = & \prefix{j}\beta \prefix{j}\overline{W}^{\dagger}  \prefix{j}\lambda 
\end{eqnarray}

Let is consider the matrix 

\beq
\prefix{j}{\cal N}^{-1} =
\prefix{j}\overline{W}^{\dagger} 
 \left[ A - \prefix{j-1}\lambda \left(E - \qn \right) \right]^{-1}
\prefix{j}\overline{V} =  \cf^{\dagger}E \prefix{j}\overline{V} =
\prefix{j}\overline{W}^{\dagger}E\ce
\eeq

It is fullfilled, because of (\ref{deduA}):

\begin{eqnarray}
\prefix{j}{\cal N}^{-1} \prefix{j}\alpha & = &
\cf^{\dagger}E \prefix{j}\overline{V} \prefix{j}\alpha \\
& = &
\prefix{j}\lambda^{-1} \cf^{\dagger}E\prefix{j}v \\
& = &
\prefix{j}\lambda^{-1} \prefix{j}\alpha
\end{eqnarray}

Analogously,

\beq
\prefix{j}\beta^{\dagger} \prefix{j}{\cal N}^{-1} =
\prefix{j}\beta^{\dagger} \prefix{j}\lambda^{-1} 
\eeq

Therefore, the matrix $\prefix{j}{\cal N}$ contains
the sought eigenvalue in its spectrum. On the other hand,
it can be shown (see appendix \ref{apeh}) that

\beq
\prefix{j}{\cal N} = \prefix{j}{\cal M}^{-1} + \prefix{j-1}\lambda I_n \label{fN}
\eeq

\noindent where $I_n$ is the identity matrix. Besides, from
(\ref{basica1}) and (\ref{basica2}),

\begin{eqnarray}
H(\prefix{j-1}\lambda ) \prefix{j}\alpha & = &
\left( \prefix{j}{\cal N} - A_{rr} \right)  \prefix{j}\alpha \\
\prefix{j}\beta^{\dagger} H(\prefix{j-1}\lambda) & = &
\prefix{j}\beta^{\dagger} \left( \prefix{j}{\cal N} - A_{rr} \right)  
\end{eqnarray}

In fact, a stronger result can be obtained:

\begin{theorem}
\label{teoHN} $H(\lambda) = {\cal N} - A_{rr}$
\end{theorem}

\begin{proof}
The proof is provided in the appendix \ref{apei}.
\end{proof}

Therefore, it is proposed the following generalization
of algorithm 2:

\vspace{3mm}
\begin{center}
\fbox{
\begin{minipage}{100cm}
\begin{tabbing}
{\bf Algorithm 7} \\
\\
{\bf Input:} $E,A, \ce, \cf$. \\
\\
{\bf Output:} $\lambda_k , \alpha_k , \beta_k$. \\
\\
{\bf 1.} \= Form $A_{rr} = \cf^{\dagger}A\ce$, and perform the eigenanalysis of $A_{rr}$, \\
{\bf 2.} \> Select the interesting modes ${}^0\lambda_k, {}^0\alpha_k, {}^0\beta_k$, \\
{\bf 3.} \> {\bf for} $j=1,2,3, \ldots$ until convergence, \\
\> {\bf 3.1} \= For each eigenvalue $k= 1,\ldots, K$, \\
\> \> {\bf 3.1.1.} \= Compute $\prefix{j}V_k = \left( A - \prefix{j-1}\lambda_k  E \right)^{-1}E \ce$, \\
\> \> {\bf 3.1.2.} \> Compute $\prefix{j}{\cal M}_k =  \cf^{\dagger}E \prefix{j}V_k$, \\
\> \> {\bf 3.1.3.} \> Compute $\prefix{j}{\cal N}_k = \prefix{j}{\cal M}_k^{-1} + \prefix{j-1}\lambda_k I_n$, \\
\> \> {\bf 3.1.4} \> Compute $\prefix{j} h_k = \left( \prefix{j}{\cal N}_k - A_{rr} \right) \prefix{j-1}\alpha_k$, \\
\> {\bf 3.2} \> Compute a matrix $\prefix{j}M$ which fullfills \\
 \>               \> $\prefix{j} M \left[ \prefix{j-1}\alpha_1, \ldots, \prefix{j-1}\alpha_K \right] =
                       \left[ \prefix{j}h_1, \ldots, \prefix{j}h_K \right]$ \\
\> {\bf 3.3} Perform the eigenanalysis of $A_{rr} + \prefix{j}M$, \\
\> {\bf 3.4} Select the interesting modes $\prefix{j}\lambda_k$, 
                    $\prefix{j}\alpha_k$, $\prefix{j}\beta_k$, \\
{\bf 4.} \> {\bf end}
\end{tabbing}
\end{minipage}
}
\end{center}
\vspace{3mm}

This algorithm is nothing else that the direct version of algorithm 2, and
it reduces to it if the $\ce$ and $\cf$ matrices are chosen as shown
in appendix \ref{apb}. A superlinear version of this
algorithm, by updating the $\ce$ and $\cf$ matrices, is also
proposed:

\vspace{3mm}
\begin{center}
\fbox{
\begin{minipage}{100cm}
\begin{tabbing}
{\bf Algorithm 8} \\
\\
{\bf Input:} $E,A, \prefix{0}\ce, \prefix{0}\cf$. \\
\\
{\bf Output:} $\lambda_k , \alpha_k , \beta_k$. \\
\\
{\bf 1.} \= Form $\prefix{0}A_{rr} = \prefix{0}\cf^{\dagger}A\prefix{0}\ce$, 
             and perform the eigenanalysis of $\prefix{0}A_{rr}$, \\
{\bf 2.} \> Select the interesting modes ${}^0\lambda_k, {}^0\alpha_k, {}^0\beta_k$, \\
{\bf 3.} \> {\bf for} $j=1,2,3, \ldots$ until convergence, \\
\> {\bf 3.1} \= For each eigenvalue $k= 1,\ldots, K$, \\
\> \> {\bf 3.1.1.} \= 
Compute $\prefix{j}V_k = \left( A - \prefix{j-1}\lambda_k  E \right)^{-1}E \prefix{j-1}\ce$, \\
\> \> {\bf 3.1.2.} \= 
Compute $\prefix{j}W_k^{\dagger} = 
\prefix{j-1}\cf^{\dagger} E \left( A - \prefix{j-1}\lambda_k  E \right)^{-1}$, \\
\> \> {\bf 3.1.3.} \> Compute $\prefix{j}{\cal M}_k =  \cf^{\dagger}E \prefix{j}V_k$, \\
\> \> {\bf 3.1.4.} \> Compute $\prefix{j}{\cal N}_k = \prefix{j}{\cal M}_k^{-1} + \prefix{j-1}\lambda_k I_n$, \\
\> \> {\bf 3.1.5} \> Compute $\prefix{j} h_k = \left( \prefix{j}{\cal N}_k - \prefix{j-1}A_{rr} \right) \prefix{j-1}\alpha_k$, \\
\> {\bf 3.2} \> Compute a matrix $\prefix{j}M$ which fullfills \\
 \>               \> $\prefix{j} M \left[ \prefix{j-1}\alpha_1, \ldots, \prefix{j-1}\alpha_K \right] =
                       \left[ \prefix{j}h_1, \ldots, \prefix{j}h_K \right]$ \\
\> {\bf 3.3} Perform the eigenanalysis of $\prefix{j-1}A_{rr} + \prefix{j}M$, \\
\> {\bf 3.4} Select the interesting modes $\prefix{j}\lambda_k$, 
                    $\prefix{j}\tilde{\alpha}_k$, $\prefix{j}\tilde{\beta}_k$, \\
\> {\bf 3.5}  \> For each eigenvalue $k= 1,\ldots, K$, compute \\
\> \> $\prefix{j}v _k= \prefix{j}V_k \prefix{j}\tilde{\alpha}_k$, and \\
\> \> $\prefix{j}w _k^{\dagger}=  \prefix{j}\tilde{\beta}_k^{\dagger}\prefix{j}W_k^{\dagger}$, \\
\> {\bf 3.6} \> Update $\prefix{j}\ce$ and $\prefix{j}\cf$ in such a way that \\
\> \> $\left[ \prefix{j}v_1, \ldots, \prefix{j}v_K \right] \in {\rm span}(\prefix{j}\ce)$, \\
\> \> $\left[ \prefix{j}w_1, \ldots, \prefix{j}w_K \right] \in {\rm span}(\prefix{j}\cf)$, \\
\>  {\bf 3.7} \> Compute $\prefix{j}A_{rr} = \prefix{j}\cf^{\dagger}A \prefix{j}\ce$, \\
\> {\bf 3.8} \> Compute $\prefix{j}\alpha_k$ and $\prefix{j}\beta_k$ by imposing \\
\> \> $\prefix{j}\ce \prefix{j}\alpha_k = \prefix{j}v_k$ and
 $\prefix{j}\cf \prefix{j}\beta_k = \prefix{j}w_k$, \\
{\bf 4.} \> {\bf end}
\end{tabbing}
\end{minipage}
}
\end{center}
\vspace{3mm}

Step 3.5 requires to pair each eigenvalue obtained of the eigeanalysis
of $A_{rr} + \prefix{j}M$ with the eigenvalues of the previous iterations.
A way to do that is to pair trying to maximize the scalar products
$\prefix{j}\tilde{\alpha}_{k1}^{\dagger}\prefix{j}\alpha_{k2}$.

\section{Numerical tests}

 The aim of this section is to apply Generalized SMA to two very
different problems: the computation of the natural frequencies
of a plate and the computation of the electromechanical modes
of an electric power system. The code was developed in
MATLAB language.

\subsection{Natural frequencies of a cross-shaped plate.}

It is intended to compute a natural frequency in a cross-shaped
plate with unequal arms. Mathematically, the problem to solve is

\beq
- \Delta \psi = \omega^2 \psi \;\; \psi(x,y) \in \Omega
\eeq

\noindent where $\Omega$ is the cross-shaped dominion,
$\omega$ the natural frequency and $\psi$ the sought mode.
In addition, the following boundary conditions must be
fulfilled:

\beq
\psi (x,y) = 0 \;\; \forall (x,y) \in \partial \Omega
\eeq

The previous partial differential equation is approximated
by a finite differences scheme. In order to apply algorithms 5, 6, 7 and 8,
it is decided that:

\begin{remunerate}

\item As it is intended to solve a symmetrical problem, matrices
          $\ce$ and $\cf$ are taken to be real and equal. Furthermore, the 
          hermitian operator $\dagger$ can de substituted by the
          transpose operator $T$, and the algorithms can be programmed
          in real, instead of complex, algebra.

\item In algorithms 5 and 7, matrix $\ce$ is a vector which approximates
           the sought mode. In algorithms 6 and 8, matrix $\ce$ is initially that
           same vector, and it is updated in each iteration to the last
          mode estimation $\prefix{j}v$.

\end{remunerate}

\begin{figure}
\hspace{30mm}
\includegraphics*[scale=.7]{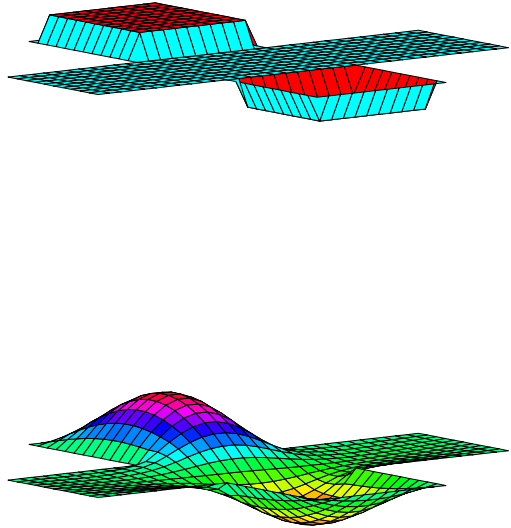}
\caption{Initial estimation and computed mode}
\label{graf}
\end{figure}

It was decided to compute the mode corresponding to the
upper arm of the cross oscillating against the lower one.
The initial mode estimation 
for algorithms 5 and 6
is shown in figure (\ref{graf}),
as well as the computed mode. 
Notice that the algorithms converge to the mode whose
shape is closest to the initial assumption, being any
initial assumption of the value of the sought eigenvalue
largely irrelevant. 

%print -dps graf.ps
% lambda = 0.11572280660335

\begin{figure}
%\hspace{-30mm}
\includegraphics*[scale=.7]{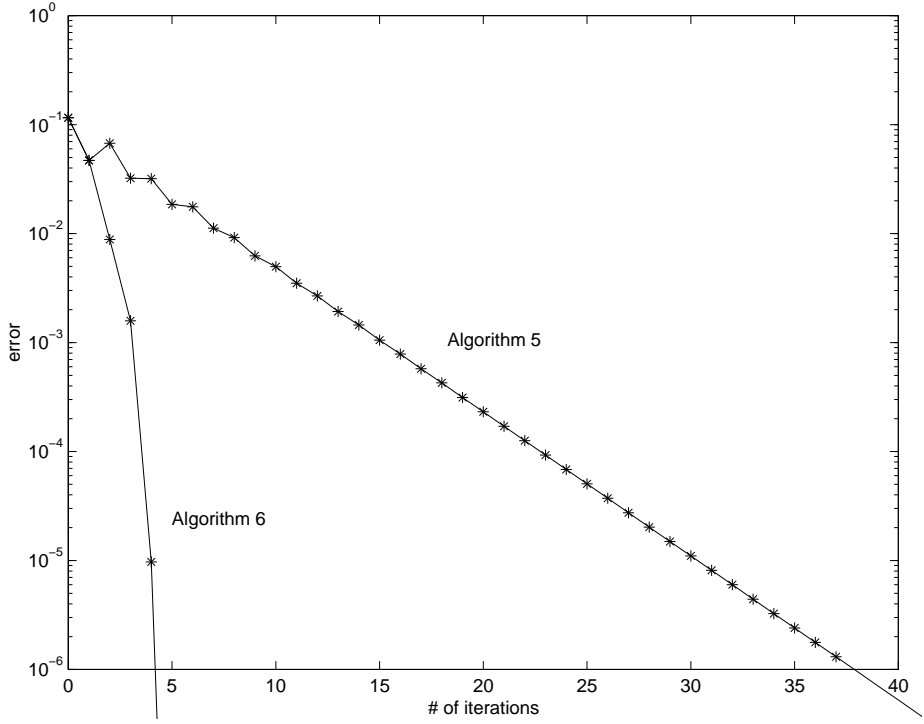}
\caption{$\omega^2$ convergence. Algorithms 5 and 6.}
\label{graf2}
\end{figure}

The eigenvalue is $\lambda = \omega^2 = 0.1157$.
Figure \ref{graf2} shows the evolution of the
absolute value of the error $\prefix{j}\lambda - \lambda$
for both algorithms 5 and 6.

\begin{figure}
%\hspace{-30mm}
\includegraphics*[scale=.7]{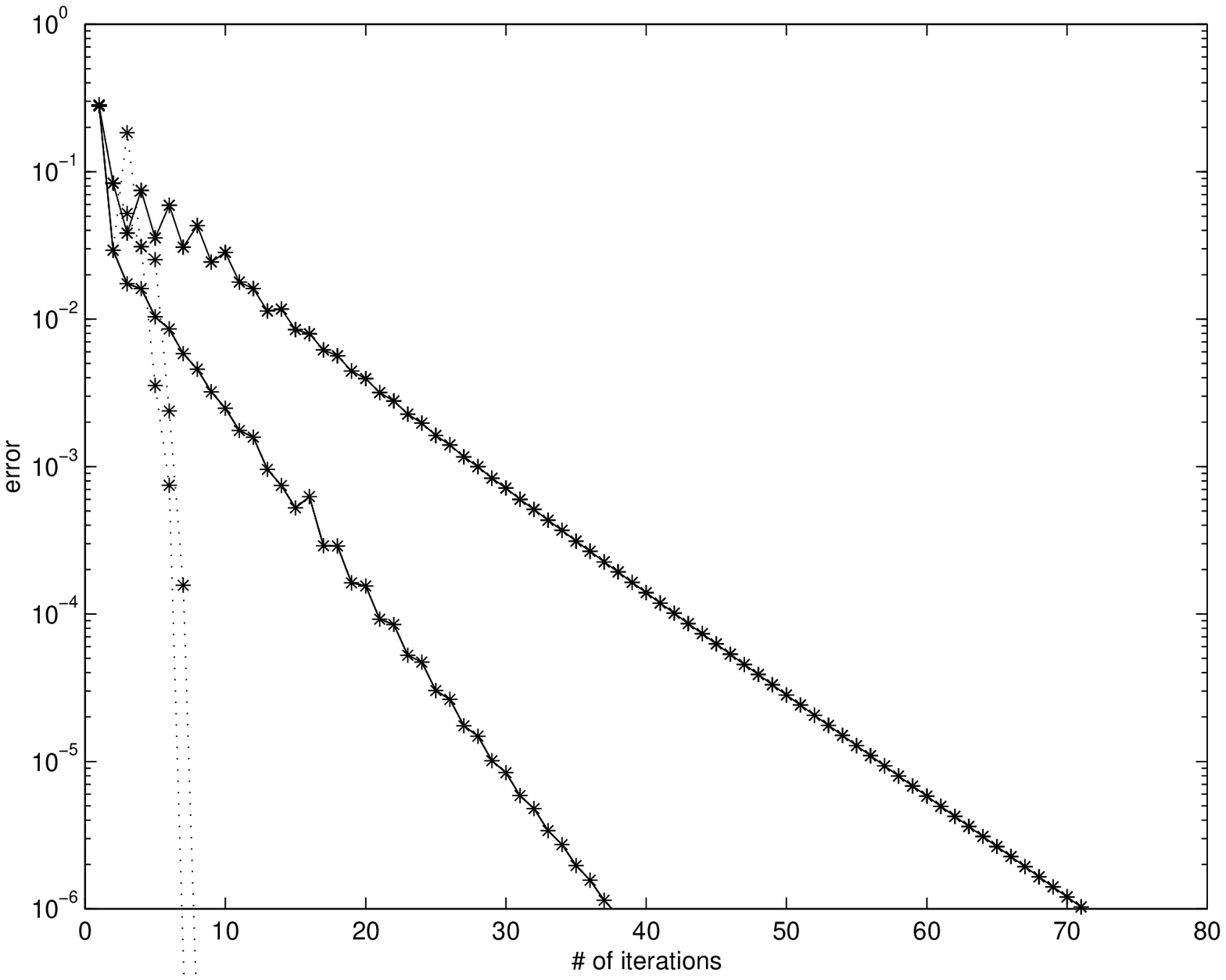}
\caption{$\omega^2$ convergence. Algorithms 7 and 8.}
\label{graf2b}
\end{figure}

To apply algorithms 7 and 8 an additional mode
$\omega^2 = 0.1243$, corresponding to the right arm
oscillating against the left arm, was computed.
Figure \ref{graf2b} shows the evolution of the
absolute value of the error (algorithm 7 in
solid lines and algorihtm 8 in dotted lines).

\subsection{Electromechanical modes
of an electric power system}

\begin{figure}
\includegraphics*[bb = 50 500 500 750, scale=.75]{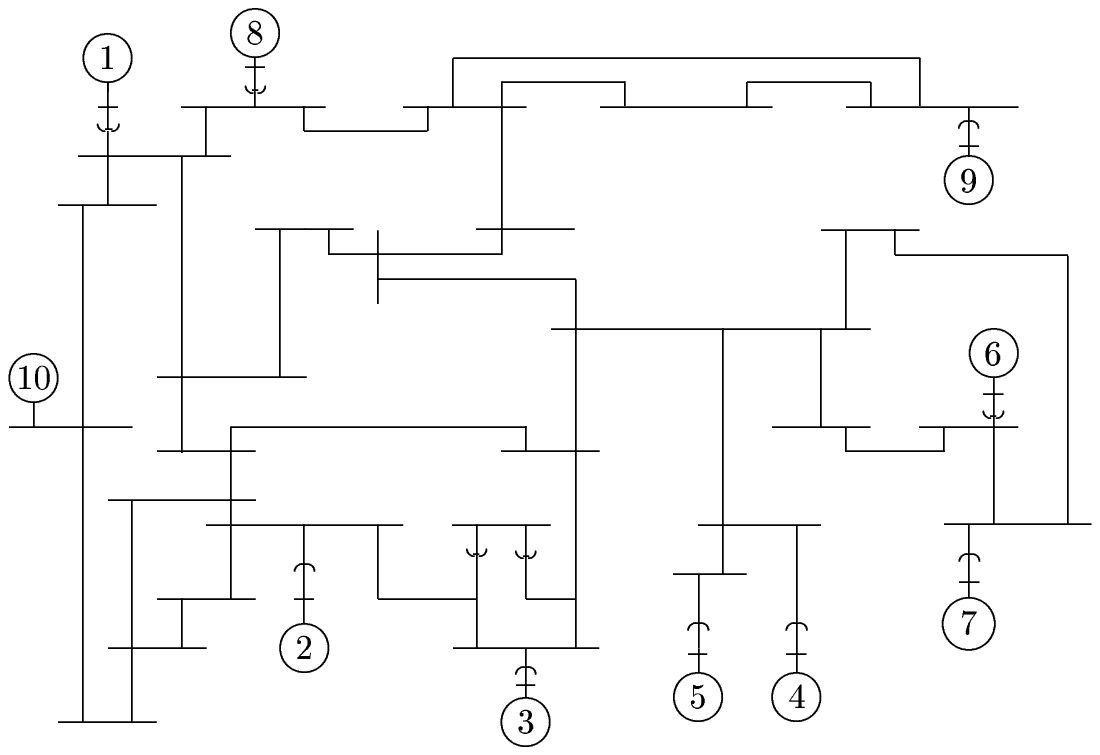}
\caption{New England Electric Power System}
\label{nered}
\end{figure}

The electric power system represented in
figure \ref{nered} (a simplified
model of the New England electric power system)
was analyzed by using SMA. The circles
reprents electric generators and the lines
the electric transmission lines. The electric generators
are modelled by 9th to 11th order linear systems,
whilst the electric network is modelled as an algebraic
constraint. Therefore, the system is a composite system,
as explained in the previous section.

The generators input variables $x_{Ik}$ are the axis and quadrature
terminal current (two variables per generator), and the output variables
$x_{Ok}$ the axis and quadrature terminal voltage (two variables
per generator). The algebraic variables $x_A$ are network voltages.
The state variables $x_{Mk}$ are mechanical, electromagnetic and
control generator variables.

The following {\em physical} information is known {\it a priori}:

\begin{remunerate}

\item The most troublesome modes are those called electromechanical
           modes. These modes are related to the generators angle $\delta_k$
           and speed $\omega_k$. It is fullfilled that $\dot{\delta_k} = 120\pi \omega_k$.
          The relationship of $\dot{\omega_k}$ with the rest of the
          variables is much more complex. However, in a very rough approximation,
          the subsystem ($\delta_k, \omega_k$) can be considered a damped
          pendulum.

\item It is known that the frequencies of the electromechanical eigenvalues
          are in the order of 1 Hz.

\item The electromechanical modes can be understood as oscillations of
           one generator or group of generators against other generator
           or group of generators.

\end{remunerate}

Fom fact 1, it follows that the $\delta_k$ and $\omega_k$ right eigencvector
components fulffill:

\beq
\lambda v(\delta_k) = 120 \pi v (\omega_k)
\eeq

On the other hand, if the rough simple pendulum model is
assumed, it should be fullfiled

\beq
\lambda^{\dagger} w(\omega_k) = 120 \pi w(\delta_k)
\eeq

Besides, it is known, from fact 2,
 that $\lambda \approx  2 \pi \imath $.
So, it is decided that, when applying algorithm 3:

\begin{eqnarray}
\ce_{Mk} & = & \left[ \begin{array}{c}
                            120 \pi \\ 2 \pi \imath \\ 0 \\ \vdots \\ 0
                            \end{array} \right] \\
 \cf_{Mk} & = & \left[ \begin{array}{c}
                            -2 \pi  \imath \\ 120 \pi  \\ 0 \\ \vdots \\ 0
                            \end{array} \right] 
\end{eqnarray}

\noindent but a normalization constant. When applying
algorithm 4, these ones are the initial values of
$\ce_{Mk}$ and $\cf_{Mk}$.

It is desired to compute the electromechanical
mode corresponding to the generators
1, 2,3, 8 and 10 oscillating againts the 4, 5, 6,
7 and 9 (the East side against the West side).
As the reducid
matrix $A_r + H(\lambda)$ is a $10 \times 10$
matrix (because each $\ce_{Mk}$ and $\cf_{Mk}$ is
a vector and there are 10 generators), it is needed to
select the relevant eigenvalue and eigenvectors
resulting from the factorization of the reduced matrix
(see Algorithm 3 and 4). The choosen procedure is
as follows:

It is defined an ``objective'' $\alpha_o$:

\beq
\alpha_o = \left[
1, 1, 1, -1, -1, -1, -1, 1, -1, 1 \right]^T
\eeq

\noindent which represents the generators oscillating
as described above. After performing the 
$\prefix{j-1}A_r + H(\prefix{j-1}\lambda)$ eigenanalysis it is obtained
the matrix spectrum $\prefix{j}\lambda_k$ and right eigenvectors
$\prefix{j}\alpha_k$, where $k=1, \ldots, 10$ (the generator number).
The products 

\beq
p_k = \alpha_o^{\dagger}\;\prefix{0}\ce_{M}^{\dagger}\;\prefix{j-1}\ce_M\;\prefix{j}\alpha_k
\eeq

\noindent are computed, and it is selected the mode which maximizes
$\| p_k \|$:

\beq
\begin{array}{c}
k = {\rm arg \;\;}{\max} \| p_k \| \\
\prefix{j}\lambda = \prefix{j}\lambda_k \\
\prefix{j}\alpha = \prefix{j}\alpha_k \\
\prefix{j}\beta = \prefix{j}\beta_k
\end{array}
\eeq

Note that, in the case of algorithm 3, the matrices $\prefix{j}\ce$ are
constant and equal to the initial one $\prefix{0}\ce$.

In algorithm 4, the vector $\prefix{j}v_M$ 
is computed in each iteration by using 
(\ref{eigc3}, \ref{decoc3b}, \ref{casf1}, \ref{casf2}),
and also $\prefix{j}v_M$ by using analogous formulae.
The matrices $\ce$ and $\cf$ are unpdated according to:

\begin{eqnarray}
\prefix{j}\ce_{Mk} & = & \prefix{j}v_{Mk} \\
\prefix{j}\cf_{Mk} & = & \prefix{j}w_{Mk} 
\end{eqnarray}

\noindent but, possibly, a normalization constant.

\begin{figure}
\hspace{10mm}
\includegraphics*[scale=.6]{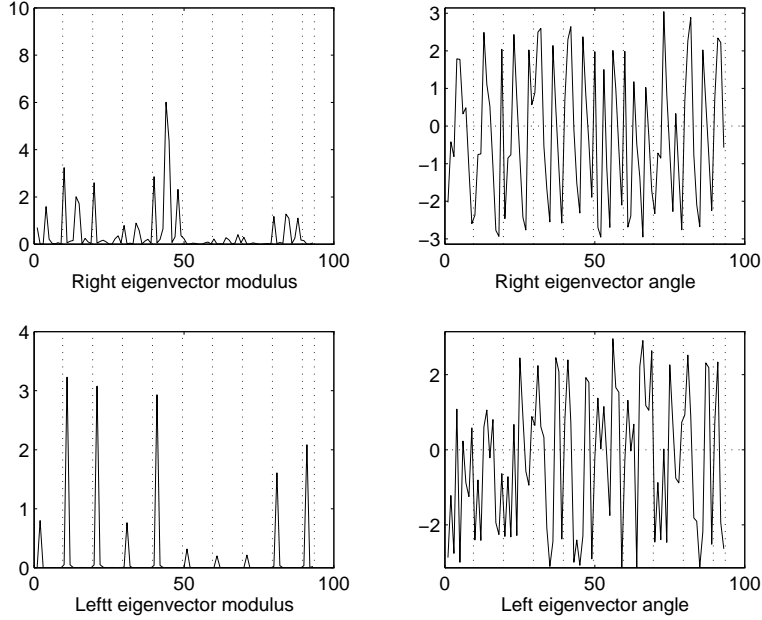}
\caption{Right and left eigenvectors}
\label{ne1}
\end{figure}

\begin{figure}
\hspace{10mm}
\includegraphics*[scale=.6]{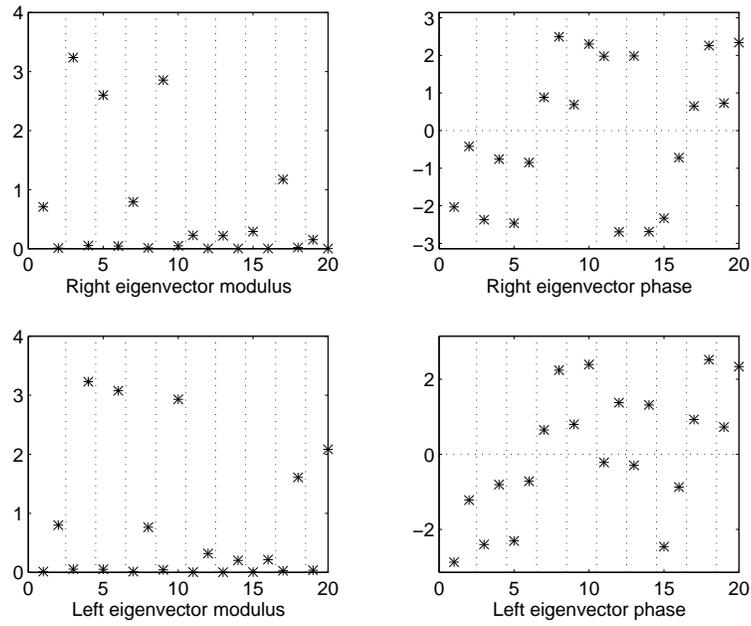}
\caption{Right and left eigenvectors. $\delta$ and $\omega$ components.}
\label{ne2}
\end{figure}

\begin{figure}
%\hspace{30mm}
\includegraphics*[scale=.7]{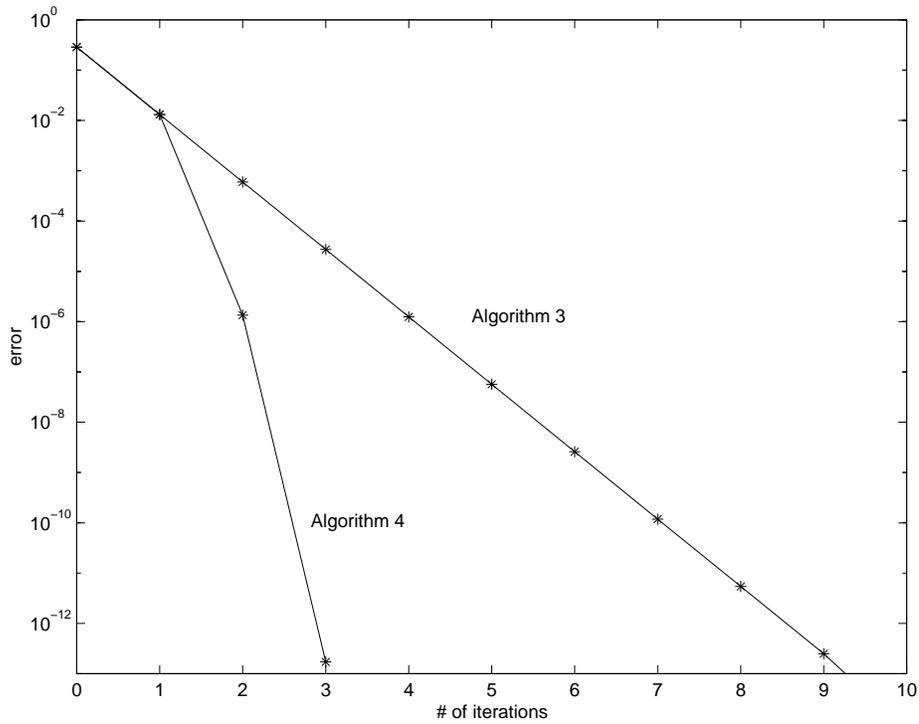}
\caption{Error evolution.}
\label{neerr}
\end{figure}

Figure \ref{ne1} shows the modulus and phase
of the right and left eigenvectors, and
figure \ref{ne2} the modulus and phase 
of the $\delta_k$ and $\omega_k$ components
of the eigenvectors. The vertical dotted lines
seggregates the variables belonging to different
generators.
Figure \ref{neerr} shows the error evolution. The computed eigenvalue
is $-0.2617 + 6.4017 \imath $.

%print -dps graf.ps
% lambda = 0.11572280660335

\section{Conclusions}

In this article  a new approach for solving
the generalized
eigenvalue problema has been introduced.
The introduced
algorithms can make efficient use
of physical information regarding the
shape of the sought eigenvectors.

\appendix
%\clearpage
\section{ Proof of the main results.}

From (\ref{eig1}) and (\ref{deco1})
it is obtained

\beq
\lambda E \ce \alpha + \lambda E z =
A \ce \alpha + A z
\label{1}
\eeq

Premultiplying (\ref{1}) by $\cf ^{\dagger}$

\beq
\lambda \cf ^{\dagger} E \ce \alpha + \lambda \cf^{\dagger} E z =
\cf ^{\dagger} A \ce \alpha + \cf ^{\dagger} A z
\eeq

But $\cf^{\dagger} E z = 0$ and $\cf^{\dagger}E\ce = I_n$. So

\beq
\lambda \alpha = \cf^{\dagger}A\ce \alpha + \cf^{\dagger}A z
\label{R}
\eeq

Note that, but the last term $\cf^{\dagger}A z$, this equation
is an eigensystem of the $n \times n$ matrix
$\cf^{\dagger}A\ce$. This matrix is usually
much smaller than $A$.

Introducing the matrices defined in (\ref{defp}):

\beq
\pn =  I_m - E \ce \cf^{\dagger} E = I_m - \qn
\eeq

These matrices are idempotent ones ($\pn = \pn^2, \qn = \qn^2$).
Premultiplying (\ref{1}) by
$\pn$, and taking into account that

\begin{eqnarray}
\pn E \ce & = & 0 \\
\pn E z & = & Ez
\end{eqnarray}

it yields

\beq
\lambda E z = \pn A \ce \alpha + \pn  A z \label{A7}
\eeq

So,

\beq
(\lambda E - A + \qn A) z = \pn A \ce \alpha
\label{bas1}
\eeq

From this equation, it is possible
to solve $z$ in function of $\alpha$.
In the same way, it is obtained

\beq
y^{\dagger} (\lambda E - A + A \qn)  = \beta^{\dagger} \cf^{\dagger} A \pn
\label{bas2}
\eeq

Taking into account that $\qn z = 0$ and $y^{\dagger}\qn = 0$, these
formulae can be written in may different ways, The most
symmetrical one is:

\begin{eqnarray}
(\lambda E - A + \left[ \qn, A \right]_+) z & = & \pn A \ce \alpha \label{bas3} \\
y^{\dagger} (\lambda E - A + \left[ \qn, A \right]_+)  & = & \beta^{\dagger} \cf^{\dagger} A \pn
\end{eqnarray}

Therefore, from (\ref{bas3}):

\begin{eqnarray}
z & = & 
\left\{ \lambda E -
A + \left[ A, \qn \right]_+ 
\right\}^{-1}\pn A\ce \alpha \label{A10} \\
y^{\dagger} & = &
\beta^{\dagger} \cf^{\dagger}  A \pn
\left\{ \lambda E -
A + \left[ A, \qn \right]_+ 
\right\}^{-1} 
\end{eqnarray}

\noindent assuming that the inverse matrix exists.
Sufficient conditions will be provided in appendix \ref{ape}.
On the other hand, it is easy to check that 
$\pn z = z, y^{\dagger} \pn = y^{\dagger}$. Therefore

\begin{eqnarray}
z & = & 
\pn\left\{ \lambda E -
A + \left[ A, \qn \right]_+
\right\}^{-1}\pn A\ce \alpha \label{zcal} \\
y^{\dagger} & = &
\beta^{\dagger} \cf^{\dagger}  A \pn
\left\{ \lambda E -
A + \left[ A, \qn \right]_+ 
\right\}^{-1} \pn \label{ycal}
\end{eqnarray}

\noindent and equation (\ref{R}) becomes

\begin{eqnarray}
\lambda \alpha & = &
\cf^{\dagger} A \ce \alpha +
\cf^{\dagger} A
\pn\left\{ \lambda E -
A + \left[ A, \qn \right]_+  
\right\}^{-1}\pn A\ce \alpha \\
& = &
A_{rr}\alpha + H(\lambda) \alpha
\end{eqnarray}

Note also, that by solving $z$ from (\ref{bas1}), it is obtained:

\beq
H(\lambda) =
\cf^{\dagger} A
\pn\left\{ \lambda E -
A + \qn A
\right\}^{-1}\pn A\ce 
\label{forref}
\eeq

\noindent and by solving $y$ from (\ref{bas2}):

\beq
H(\lambda) =
\cf^{\dagger} A
\pn\left\{ \lambda E -
A + A \qn
\right\}^{-1}\pn A\ce 
\eeq

These ones are just some few of the many equivalent
ways to write $H(\lambda)$. Some of then many be more
amenable for computation than (\ref{T2}).

 \section{The relationship with ``classical'' SMA}
\label{apb}

 The formulae (\ref{T1},\ref{T2}) are just the
generalized version of the ``clasical''
SMA (Selective Modal
Analysis) 
formulae, as defined in section 1.

 To show the relationship, consider the 
problem (\ref{uno}). Let also assume that $A$,
$v$, $w$ are partitioned according
(\ref{tres},\ref{tresb}). Therefore, the
relevant subspaces in the generalized version are:

\beq
\ce = \left[ \begin{array}{cccc}
             1 & 0 & \ldots & 0 \\
             0 & 1 & \ldots & 0 \\
           \vdots & \vdots & \ddots & \vdots \\
             0 & 0 & \ldots & 1 \\
             0 & 0 & \ldots & 0 \\
           \vdots & \vdots & \ddots & \vdots \\
             0 & 0 & \ldots & 0 
             \end{array} \right] 
=
\left[ \begin{array}{c}
I_n \\ 0 
\end{array} \right] 
=  \cf
\eeq

So

\beq
\cf^{\dagger} A \ce = 
\left[ \begin{array}{cc}
I_n & 0 
\end{array} \right] 
\left[ \begin{array}{cc}
            A_{rr} & A_{rz} \\ A_{zr} & A_{zz}
            \end{array} \right]
\left[ \begin{array}{c}
I_n \\ 0 
\end{array} \right] 
= A_{rr}
\eeq

\noindent as suggested by the notation. On the other hand

\beq
\ce \cf^{\dagger} =
\qn =
\left[ \begin{array}{c}
I_n \\ 0 
\end{array} \right] 
\left[ \begin{array}{cc}
I_n & 0 
\end{array} \right] 
=
\left[ \begin{array}{cc}
I_n & 0 \\ 0 & 0
\end{array} \right] 
\eeq

Therefore

\begin{eqnarray}
& &
\left[ A, \qn \right]_+ \nonumber \\
& = & 
\left[ \begin{array}{cc}
            A_{rr} & A_{rz} \\ A_{zr} & A_{zz}
            \end{array} \right]
\left[ \begin{array}{cc}
I_n & 0 \\ 0 & 0
\end{array} \right] 
+
\left[ \begin{array}{cc}
I_n & 0 \\ 0 & 0
\end{array} \right] 
\left[ \begin{array}{cc}
            A_{rr} & A_{rz} \\ A_{zr} & A_{zz}
            \end{array} \right] \nonumber \\
& = &
\left[ \begin{array}{cc}
            A_{rr} & 0 \\ A_{zr} & 0
            \end{array} \right]
+
\left[ \begin{array}{cc}
            A_{rr} & A_{rz} \\ 0 & 0
            \end{array} \right]
\nonumber \\
& = &
\left[ \begin{array}{cc}
            2A_{rr} & A_{rz} \\ A_{zr} & 0
            \end{array} \right]
\end{eqnarray}

Then,

\begin{eqnarray}
\lambda -  A + \left[ A, \qn \right]_+
 & = & 
\left[ \begin{array}{cc}
            \lambda + A_{rr} & 0 \\ 0 & \lambda - A_{zz}
            \end{array} \right] \\
\left\{
\lambda -  A + \left[ A, \qn \right]_+
\right\} ^{-1} & = & 
\left[ \begin{array}{cc}
                         \left( \lambda + A_{rr} \right)^{-1}
            & 0 \\ 0 & 
            \left( \lambda - A_{zz} \right)^{-1}
            \end{array} \right] 
\end{eqnarray}

Besides

\beq
\pn = I_m - \qn =
\left[ \begin{array}{cc}
            0 & 0 \\ 0 & I_{m-n}
            \end{array} \right] 
\eeq

So

\begin{eqnarray}
\pn A \ce & = & 
\left[ \begin{array}{cc}
            0 & 0 \\ 0 & I_{m-n}
            \end{array} \right] 
\left[ \begin{array}{cc}
            A_{rr} & A_{rz} \\ A_{zr} & A_{zz}
            \end{array} \right] 
\left[ \begin{array}{c}
I_n \\ 0 
\end{array} \right] \nonumber \\
& = &
\left[ \begin{array}{cc}
            0 & 0 \\ 0 & I_{m-n}
            \end{array} \right] 
\left[ \begin{array}{c}
A_{rr} \\ A_{zr} 
\end{array} \right] \nonumber \\
& = &
\left[ \begin{array}{c}
0 \\ A_{zr} 
\end{array} \right] \\
\cf^{\dagger}A \pn & = &
\left[ \begin{array}{cc}
I_n & 0 
\end{array} \right] 
\left[ \begin{array}{cc}
            A_{rr} & A_{rz} \\ A_{zr} & A_{zz}
            \end{array} \right] 
\left[ \begin{array}{cc}
            0 & 0 \\ 0 & I_{m-n}
            \end{array} \right] \nonumber \\
& = &
\left[ \begin{array}{cc}
A_{rr} & A_{rz} 
\end{array} \right] 
\left[ \begin{array}{cc}
            0 & 0 \\ 0 & I_{m-n}
            \end{array} \right] \nonumber \\
& = &
\left[ \begin{array}{cc}
0 & A_{rz} 
\end{array} \right] 
\end{eqnarray}

And, finally

\begin{eqnarray}
H(\lambda) & = &
\left( \cf^{\dagger} E \tilde{\pn} \right)
\left\{
\lambda -  A + \left[ A, \qn \right]_+
 \right\} ^{-1} 
\left( \tilde{\pn} A \ce \right) \nonumber \\
& = &
\left[ \begin{array}{cc}
0 & A_{rz} 
\end{array} \right] 
\left[ \begin{array}{cc}
                        \left( \lambda + A_{rr} \right)^{-1}
 & 0 \\ 0 & 
            \left( \lambda - A_{zz} \right)^{-1}
            \end{array} \right] 
\left[ \begin{array}{c}
0 \\ A_{zr} 
\end{array} \right] \nonumber \\
& = &
A_{rz} \left( \lambda - A_{zz} \right)^{-1} A_{zr}
\end{eqnarray}

\section{Proof of theorem \ref{teo1}}

\begin{proof}
Let us recall $\prefix{j}\epsilon = \prefix{j}\lambda - \lambda$.
$\lambda$ is in the spectrum of $A_{rr} + H(\lambda)$
and, because algorithm 3, $\prefix{j}\lambda$ is in
the spectrum of $A_{rr} + H(\prefix{j-1}\lambda)$. 
But

\beq
A_{rr} + H(\prefix{j-1}\lambda) =
A_{rr} + H(\lambda) + \frac{\partial H(\lambda)}{\partial \lambda}\; \prefix{j-1}\epsilon
+ o(\prefix{j-1}\epsilon)
\eeq

By appliying a well-known perturbation
formula:

\beq
\prefix{j}\lambda - \lambda = \prefix{j}\epsilon =
\frac{\beta^{\dagger}\frac{\partial H(\lambda)}{\partial \lambda}\alpha}{\beta^{\dagger} \alpha}
\;\prefix{j-1}\epsilon + o(\prefix{j-1}\epsilon) \label{lineal}
\eeq

Let us define

\beq
\rho^{-1} = -\frac{\beta^{\dagger}\frac{\partial H(\lambda)}{\partial \lambda}\alpha}{\beta^{\dagger} \alpha}
\eeq

If $\lambda$ and $\prefix{j-1}\lambda$ are close enough, it
is possible to neglect the higher order terms $o(\epsilon^{j-1})$.
Then, it is clear that the algorithm converges 
($\epsilon^j \rightarrow 0$) if and only if $\mid \rho^{-1} \mid < 1$.

On the other hand, from
(\ref{T2}), (\ref{zcal}) and (\ref{ycal}):

\begin{eqnarray}
-\beta^{\dagger}\frac{\partial H(\lambda)}{\partial \lambda}\alpha
& = & \nonumber \\
\beta
\cf^{\dagger} A \pn \left\{ \lambda E - A + 
\left[ A, \qn \right]_+
\right\}^{-1} E
\left\{ \lambda E - A + 
\left[ A, \qn \right]_+
\right\}^{-1}
\pn A \ce \alpha 
& = & \nonumber \\
y^{\dagger}Ez & &
\end{eqnarray}

Therefore

\beq
\rho = \frac{\beta^{\dagger} \alpha}{y^{\dagger} E z}
\eeq

But, premultiplying (\ref{deco1}) and
posmultiplying (\ref{deco3}) by $\qn$:

\begin{eqnarray}
\qn v & = & \qn \ce \alpha + \qn z \nonumber  \\
& = & E \ce (\cf ^{\dagger} E \ce) \alpha + E \ce (\cf ^{\dagger} E z) \nonumber \\
& = & E \ce \alpha \\
w^{\dagger} \qn & = & \beta^{\dagger} \cf^{\dagger} + y^{\dagger} \nonumber \\
& = & \beta^{\dagger} (\cf ^{\dagger} E \ce) \cf^{\dagger} E + (y^{\dagger} E \ce) \cf^{\dagger} E \nonumber \\
& = & \beta^{\dagger} \cf^{\dagger} E
\end{eqnarray}

So

\beq
w^{\dagger} \qn v = w^{\dagger}\qn \qn v = 
\beta^{\dagger} \cf^{\dagger} E E \ce \alpha = 
\beta^{\dagger} \alpha
\eeq

On the other hand,

\begin{eqnarray}
w^{\dagger} E v & = & (\beta^{\dagger} \cf^{\dagger} + y^{\dagger})E(\ce \alpha + z) \nonumber \\
& = & \beta^{\dagger} \cf^{\dagger} E \ce \alpha +
         \beta^{\dagger} ( \cf^{\dagger} E z) + (y^{\dagger} E \ce) \alpha +
         y^{\dagger} E z \nonumber \\
& = & \beta^{\dagger} \alpha + y^{\dagger} E z \nonumber \\
& = & w^{\dagger} \qn v + y^{\dagger} E z
\end{eqnarray}

So

\beq
y^{\dagger} E z = w^{\dagger} (E - \qn) v
\eeq

\noindent and

\beq
\rho = \frac{w^{\dagger} \qn v}{w^{\dagger} (E- \qn) v}
\label{C18}
\eeq

\end{proof}

\section{Proof of theorem \ref{teoce}}
\label{apecc}

\begin{proof}
The invariance with respect to transformations
$\ce \leftarrow \ce + (I_m - E) \cl$ shall be proven,
being the other case esentially equal.

Firstly, note that $\qn$ and $\pn$ are invariant
under the transformation

\beq
E\left(\ce + (I_m-E) \cl\right)\cf^{\dagger}E =
E\ce\cf^{\dagger}E + E (I_m - E)\cl\cf^{\dagger}E =
E\ce\cf^{\dagger}E 
\eeq

\noindent which proves the $\qn$ invariance.
As $\pn = I_m - \qn$, $\pn$ is also invariant.
Let us consider now:

\begin{eqnarray}
\pn A \left(I_m - E\right) \cl & = & \pn A \left(I_m - E\right) \cl  
\label{des1} \\
-\left\{ -A + \qn A \right\}
\left(I_m - E\right) \cl & = & \pn A \left(I_m - E\right) \cl  
\label{des2} \\
-\left\{\lambda E -A + \qn A \right\}
\left(I_m - E\right) \cl & = & \pn A \left(I_m - E\right) \cl  
\label{des3} \\
-\left(I_m - E\right) \cl & = & 
\left\{\lambda E -A + \qn A \right\}^{-1}\pn A \left(I_m - E\right) \cl  
\label{des4} \\
-\left(I_m - E\right) \cl & = & 
\pn \left\{\lambda E -A + \qn A \right\}^{-1}\pn A \left(I_m - E\right) \cl  
\label{des5} \\
-\cf^{\dagger}A\left(I_m - E\right) \cl & = & 
\cf^{\dagger}A
\pn \left\{\lambda E -A + \qn A \right\}^{-1}\pn A \left(I_m - E\right) \cl  
\nonumber \\
\label{des6} 
\end{eqnarray}

It has been used that $\pn = I_m - \qn$ to go from (\ref{des1})
to (\ref{des2}), 
that $E(I_m - E) = 0$ to go from (\ref{des2}) to (\ref{des3})
and that
that $\pn (I_m - E) = (I_m - \qn)(I_m - E) = I_m - E$
to go from (\ref{des4}) to (\ref{des5}).

On the other hand under the 
considered transformation
$A_{rr} + H(\lambda)$ becomes

\begin{eqnarray}
\cf^{\dagger}A \left( \ce + (I_m - E)\cl \right)
+ 
\cf^{\dagger} A \pn 
\left\{\lambda E -A + \qn A \right\}^{-1}\pn A 
\left( \ce + (I_m - E)\cl \right) & = & \nonumber \\
\cf^{\dagger}A \ce
+ 
\cf^{\dagger} A \pn 
\left\{\lambda E -A + \qn A \right\}^{-1}\pn A 
\ce  &  & \label{des7} \\
+ \left[
\cf^{\dagger}A (I_m - E)\cl 
+
\cf^{\dagger} A \pn 
\left\{\lambda E -A + \qn A \right\}^{-1}\pn A 
 (I_m - E)\cl \right] & = & \nonumber \\
\cf^{\dagger}A \ce
+ 
\cf^{\dagger} A \pn 
\left\{\lambda E -A + \qn A \right\}^{-1}\pn A 
\ce  &  & \label{des8} 
\end{eqnarray}

It has been used the $H(\lambda)$ formula
(\ref{forref}) to write (\ref{des7}), and
(\ref{des6}) to simplify it.

\end{proof}

\section{Proof of theorem \ref{teo2}}
\label{aped}

The net effect of steps 3.4 and 3.5
is to impose

\beq
\prefix{j}\cn_v
\prefix{j+1}v = \left( I + \prefix{j}\pn N(\prefix{j}\lambda, \prefix{j}\qn)
                                          \prefix{j}\pn A \right) \prefix{j}v
\eeq

This is because matrix $\prefix{j}\ce$ is proportional to the vector
$\prefix{j}v$. $\prefix{j}\cn_v$ is a normalization constant. In fact, algorithm
4 is invariant with respect to arbitrary normalizations of vectors
$\prefix{j}v$ and $\prefix{j}w$ (only the normalized matrices 
$\prefix{j}\ce$, $\prefix{j}\cf$, $\prefix{j}\qn$ and $\prefix{j}\pn$ are
required). 

In order to fix the normalization $\prefix{j}\cn_v$, let us focus in equations
(\ref{deco1}) and (\ref{zcal}).

\beq
v = \left( I + \prefix{j}\pn N(\lambda, \prefix{j}\qn)
                                          \prefix{j}\pn A \right) \prefix{j}v
\label{vdef}
\eeq

$v$ is a right eigenvector. Then, the normalization constant
$\prefix{j}\cn_v$ is implicitly chosen by enacting that this formula
is valid for all $i$, with the same eigenvector (i.e., with the
same phase and absolute value). Of course. there are analogous
formulae for the left eigenvectors. Furthermore, it is required
that

\beq
w^{\dagger}E v = 1
\eeq

Now, from (\ref{des4}) it can be
deduced that

\begin{eqnarray}
\prefix{j}\cn_v
E \prefix{j+1}v & = & \left( E + E \prefix{j}\pn N(\prefix{j}\lambda, \prefix{j}\qn)
                                          \prefix{j}\pn A \right) E \prefix{j}v \label{vdd1} \\
Ev  & = & \left( E + E \prefix{j}\pn N(\lambda, \prefix{j}\qn)
                                          \prefix{j}\pn A \right) E \prefix{j}v \label{vdd2} \\
& = & E\prefix{j}v + E\prefix{j}z
\end{eqnarray}

The last line uses the error $\prefix{j}z = v - \prefix{j}z$,

\begin{eqnarray}
\prefix{j}z & = & \prefix{j}\pn N(\lambda, \prefix{j}\qn)
                                          \prefix{j}\pn A  \prefix{j}v \label {prim1} \\
& = & N(\lambda, \prefix{j}\qn)
                                          \prefix{j}\pn A   \prefix{j}v  \label{prim2}
\end{eqnarray}

The last equation follows from (\ref{A10}).
Let us define the eigenvalue error $\prefix{j}\epsilon =
\lambda - \prefix{j}\lambda$. It is easy
to check that

\begin{eqnarray}
\frac{\partial N(\lambda,\prefix{j}\qn) }{\partial \lambda} & = & - N(\lambda,\prefix{j}\qn) E N(\lambda,\prefix{j}\qn)   
\label{prder} \\
\frac{\partial^2 N(\lambda,\prefix{j}\qn) }{\partial \lambda^2} & = & 2 N(\lambda,\prefix{j}\qn)  E N(\lambda,\prefix{j}\qn) 
 E N(\lambda,\prefix{j}\qn)  \\
\vdots & = & \vdots \nonumber
\end{eqnarray}

Then
substracting (\ref{vdd2}) from (\ref{vdd1}):

\begin{eqnarray}
\hspace{4mm} & &  
\prefix{j}\cn_v E \prefix{j+1}v - Ev \nonumber \\
& = &
E \left\{
\prefix{j}\pn N(\prefix{j}\lambda, \prefix{j}\qn) -
\prefix{j}\pn N(\lambda, \prefix{j}\qn)
\right\} E \prefix{j}\pn A \prefix{j}v  \nonumber \\
& = &
E \left\{
-\prefix{j}\pn N(\lambda, \prefix{j}\qn) E N(\lambda, \prefix{j}\qn)
\prefix{j}\epsilon + \right.   \nonumber \\ 
& & \left.
\prefix{j}\pn N(\lambda, \prefix{j}\qn) E N(\lambda, \prefix{j}\qn)
E N(\lambda, \prefix{j}\qn)
\prefix{j}\epsilon^2 + \ldots \right\}
E \prefix{j}\pn A \prefix{j}v  \nonumber \\
& = &
E \left\{
-\prefix{j}\pn N(\lambda, \prefix{j}\qn) \prefix{j}\epsilon + \right.  \nonumber \\ 
& & \left. 
\prefix{j}\pn N(\lambda, \prefix{j}\qn) E N(\lambda, \prefix{j}\qn)
\prefix{j}\epsilon^2 + \ldots \right\}
E N(\lambda, \prefix{j}\qn)
E \prefix{j}\pn A \prefix{j}v \nonumber \\
& = &
E \left\{
-\prefix{j}\pn N(\lambda, \prefix{j}\qn) \prefix{j}\epsilon +
\prefix{j}\pn N(\lambda, \prefix{j}\qn) E N(\lambda, \prefix{j}\qn)
\prefix{j}\epsilon^2 + \ldots \right\}
E \prefix{j}z  \nonumber \\
& = &
E \left\{
-\prefix{j}\pn N(\lambda, \prefix{j}\qn) \prefix{j}\epsilon +
\prefix{j}\pn N(\lambda, \prefix{j}\qn) E N(\lambda, \prefix{j}\qn)
\prefix{j}\epsilon^2 + \ldots \right\}
\prefix{j}\pn E \prefix{j}z \nonumber \\
& = &
E \left\{
-E \prefix{j}\pn N(\lambda, \prefix{j}\qn)\prefix{j}\pn E \prefix{j}\epsilon + \right. \nonumber \\
& & \left.
E \prefix{j}\pn N(\lambda, \prefix{j}\qn) E N(\lambda, \prefix{j}\qn)
\prefix{j}\pn E \prefix{j}\epsilon^2 + \ldots \right\}
\prefix{j}\pn E \prefix{j}z \nonumber \\
& = &
E R_z(\prefix{j}\qn, \prefix{j}\epsilon) E \prefix{j}z \label{C}
\end{eqnarray}

In order to derive the last three equations it has been used
that $\prefix{j}\pn \prefix{j}z = \prefix{j}z$
and that $\prefix{j}\pn E = \left( \prefix{j}\pn E \right)^2$. 
The last equation is the definition of matrix $R_z$. 
This matrix is well defined so long as the matrix
$N(\lambda, \prefix{j}\qn)$ is bounded and 
$\prefix{j}\epsilon$ is small enough. Sufficient conditions
shall be discussed later on. On the other hand,
it is obvious that:

\beq
R_z(\prefix{j}\qn, \prefix{j}\epsilon)  =
\prefix{j}\pn R_z(\prefix{j}\qn, \prefix{j}\epsilon)  =
R_z(\prefix{j}\qn, \prefix{j}\epsilon) \prefix{j}\pn
\eeq

Let us denote by $\qn$ the matrix formed
with the eigenvectors

\beq
\qn = Ev w^{\dagger}E
\eeq

Then 

\begin{eqnarray}
E R_z(\prefix{j}\qn, \prefix{j}\epsilon)  \qn \prefix{j}z & = &
 E R_z(\prefix{j}\qn, \prefix{j}\epsilon) \prefix{j}\pn  \qn \prefix{j}z \\
& = &
 E R_z(\prefix{j}\qn, \prefix{j}\epsilon) \prefix{j}\pn  E v w^{\dagger}E  \prefix{j}z \nonumber \\
& = &
\left(  w^{\dagger}E  \prefix{j}z \right)
 E R_z(\prefix{j}\qn, \prefix{j}\epsilon) \prefix{j}\pn  E \left( \prefix{j}v + \prefix{j}z \right) \nonumber \\
& = &
\left(  w^{\dagger}E  \prefix{j}z \right)
 E R_z(\prefix{j}\qn, \prefix{j}\epsilon) \prefix{j}\pn  E \prefix{j}z  \nonumber \\
& = &
\left(  w^{\dagger}E  \prefix{j}z \right)
 E R_z(\prefix{j}\qn, \prefix{j}\epsilon) E \prefix{j}z  \nonumber
\end{eqnarray}

So

\beq
E R_z(\prefix{j}\qn, \prefix{j}\epsilon) \left( E - \qn \right) \prefix{j}z =
\left[ 1 -
\left(  w^{\dagger}E  \prefix{j}z \right) \right]
 E R_z(\prefix{j}\qn, \prefix{j}\epsilon) E \prefix{j}z  
\label{jba}
\eeq

On the other hand

\begin{eqnarray}
\prefix{j}\cn_v E \prefix{j+1}v - Ev & = &
\prefix{j}\cn_v E \left( v - \prefix{j+1}z \right) - Ev \nonumber \\
& = &
\left( \prefix{j}\cn_v -1 \right) E v - 
\prefix{j}\cn_v E \prefix{j+1}z \nonumber \\ & = &
 E R_z(\prefix{j}\qn, \prefix{j}\epsilon) E \prefix{j}z  
\end{eqnarray}

So, premultiplying by $E - \qn$, taking into account that
$(E - \qn) v = 0$ and (\ref{jba}), it yields

\beq
( E - \qn) \prefix{j+1}z = 
\frac{-  (E - \qn) R_z(\prefix{j}\qn, \prefix{j}\epsilon) (E - \qn) \prefix{j}z}
        {\prefix{j}\cn_v \left[ 1 - \left(  w^{\dagger}E  \prefix{j}z \right) \right]}
\label{zpver1}
\eeq

Premultiplying (\ref{C})

\begin{eqnarray}
\prefix{j}\cn_v w^{\dagger}E \prefix{j+1}v & = &
1 + w^{\dagger}E
R_z(\prefix{j}\qn, \prefix{j}\epsilon) E \prefix{j}z \nonumber \\
& = &
1 + \frac{w^{\dagger}ER_z(\prefix{j}\qn, \prefix{j}\epsilon)
\left( E - \qn \right) \prefix{j}z}
{\left[ 1 - \left(  w^{\dagger}E  \prefix{j}z \right) \right]} \label{C1}
\end{eqnarray}

But,

\beq
w^{\dagger}E \prefix{j}\pn =
(\prefix{j}w^{\dagger} + \prefix{j}y^\dagger)
E \prefix{j}\pn = \prefix{j}y^\dagger
E \prefix{j}\pn 
\eeq

In this formula, it has been used the error
$\prefix{j}y = w - \prefix{j}w$.  The dual
equations of (\ref{prim1},\ref{prim2}) are:

\begin{eqnarray}
\prefix{j}y^{\dagger} & = & \prefix{j}w^{\dagger} A \prefix{j}\pn N(\lambda, \prefix{j}\qn)
                                          \prefix{j}\pn \label {dual1} \\
& = & \prefix{j}w^{\dagger} A \prefix{j}\pn N(\lambda, \prefix{j}\qn)
                                         \label{dual2}
\end{eqnarray}

Besides

\beq
\prefix{j}y^{\dagger} \qn \prefix{j}\pn =
( \prefix{j}y^{\dagger} E v) (w^{\dagger}E \prefix{j}\pn) =
( \prefix{j}y^{\dagger} E v) (\prefix{j}y^{\dagger}E \prefix{j}\pn) 
\eeq

So

\beq
\prefix{j}y^{\dagger} (E - \qn) \prefix{j}\pn =
\left[1 -  ( \prefix{j}y^{\dagger} E v) \right]
\prefix{j}y^{\dagger}E \prefix{j}\pn
\eeq

\noindent and

\beq
w^{\dagger}E \prefix{j}\pn
= \frac{\prefix{j}y^{\dagger} (E - \qn) \prefix{j}\pn}
{\left[1 -  ( \prefix{j}y^{\dagger} E v) \right]}
\eeq

So, as 

\begin{eqnarray}
w^{\dagger}E R_z(\prefix{j}\qn, \prefix{j}\epsilon) & = &
w^{\dagger}E \prefix{j}\pn R_z(\prefix{j}\qn, \prefix{j}\epsilon) \nonumber \\
& = &
\frac{\prefix{j}y^{\dagger} (E - \qn) \prefix{j}\pn R_z(\prefix{j}\qn, \prefix{j}\epsilon) }
{\left[1 -  ( \prefix{j}y^{\dagger} E v) \right]} \nonumber \\
& = &
\frac{\prefix{j}y^{\dagger} (E - \qn)  R_z(\prefix{j}\qn, \prefix{j}\epsilon) }
{\left[1 -  ( \prefix{j}y^{\dagger} E v) \right]} 
\end{eqnarray}

\noindent equation (\ref{C1}) yields:

\beq
\prefix{j}\cn_v w^{\dagger}E\prefix{j+1}v = 
1 +
\frac{\prefix{j}y^{\dagger} (E - \qn)  R_z(\prefix{j}\qn, \prefix{j}\epsilon) (E- \qn) \prefix{j}z}
{\left[1 -  ( \prefix{j}y^{\dagger} E v) \right] \left[ 1 - \left(  w^{\dagger}E  \prefix{j}z \right) \right]} 
\label{normv1}
\eeq

But

\beq
1 - w^{\dagger}E \prefix{j}z =
w^{\dagger}Ev -  w^{\dagger}E \prefix{j}z =
w^{\dagger}E\prefix{j}v
\eeq

\noindent and

\beq
1 - \prefix{j}y^{\dagger}Ev = 
w^{\dagger}Ev - \prefix{j}y^{\dagger}Ev = 
\prefix{j}w^{\dagger}Ev
\eeq

On the other hand, from (\ref{prim1}) and (\ref{dual1}),

\beq
\prefix{j}w^{\dagger}E \prefix{j}z = 0 = \prefix{j}y^{\dagger}E \prefix{j}v
\label{conper}
\eeq

So,

\beq
\prefix{j}w^{\dagger}Ev
=
\prefix{j}w^{\dagger}E(\prefix{j}v + \prefix{j}z) =
\prefix{j}w^{\dagger}E\prefix{j}v 
=
w^{\dagger}E\prefix{j}v
\eeq

Therefore, from (\ref{normv1})

\beq
\prefix{j}\cn_v \prefix{j+1}w^{\dagger}E\prefix{j+1}v = 
1 +
\frac{\prefix{j}y^{\dagger} (E - \qn)  R_z(\prefix{j}\qn, \prefix{j}\epsilon) (E- \qn) \prefix{j}z}
{ \left( \prefix{j}w^{\dagger}E\prefix{j}v \right)^2}
\label{normv2}
\eeq

Subtituting in (\ref{zpver1})

\beq
\frac{1}{\prefix{j+1}wE\prefix{j+1}v}
( E - \qn) \prefix{j+1}z = 
\frac{-  (E - \qn) R_z(\prefix{j}\qn, \prefix{j}\epsilon) }
         {1 + 
         \frac{\prefix{j}y^{\dagger} (E - \qn)  R_z(\prefix{j}\qn, \prefix{j}\epsilon) (E- \qn) \prefix{j}z}
        { \left( \prefix{j}w^{\dagger}E\prefix{j}v \right)^2}}
\frac{( E - \qn) \prefix{j}z  }{\prefix{j}wE\prefix{j}v}
\label{zpver2}
\eeq

Let us define the vectors

\begin{eqnarray}
\prefix{j}\tilde{z} & = & \frac{1}{\prefix{j}wE\prefix{j}v} \left( E - \qn \right) \prefix{j}z  \label{E50}\\
\prefix{j}\tilde{y}^{\dagger} & = & \frac{1}{\prefix{j}wE\prefix{j}v} \prefix{j}y^{\dagger} 
\left( E - \qn \right) \label{E51}
\end{eqnarray}

Then, equation (\ref{zpver2}) can be written as

\beq
\prefix{j+1}\tilde{z} = \frac{- (E- \qn) R_z(\prefix{j}\qn, \prefix{j}\epsilon) }
                                                { 1 + \prefix{j}\tilde{y}^{\dagger}
                                                   R_z(\prefix{j}\qn, \prefix{j}\epsilon) \prefix{j}\tilde{z}}
                                        \prefix{j}\tilde{z}
\label{zti}
\eeq

There is also a dual equation

\beq
\prefix{j+1}\tilde{y}^{\dagger} = \prefix{j}\tilde{y}^{\dagger}
\frac{-  R_z(\prefix{j}\qn, \prefix{j}\epsilon) (E -  \qn)}
                                                { 1 + \prefix{j}\tilde{y}^{\dagger}
                                                   R_z(\prefix{j}\qn, \prefix{j}\epsilon) \prefix{j}\tilde{z}}
\label{yti}
\eeq

Now, let us write the partitipation factor $\prefix{j}\rho$ in terms
of these new error vectors. From (\ref{C18})

\beq
\prefix{j}\rho^{-1} = \frac{w^{\dagger} ( E - \prefix{j}\qn ) v}{w^{\dagger} \prefix{j}\qn v}
\eeq

Taking into account that

\beq
\prefix{j}\qn = \frac{1}{\prefix{j}w^{\dagger}E\prefix{j}v}E\prefix{j}v\prefix{j}w^{\dagger}E
\eeq

\noindent it is easy to check that

\beq
\prefix{j}\rho^{-1} = \frac{\prefix{j}w^{\dagger} ( E - \qn ) \prefix{j}v}
                                             {\prefix{j}w^{\dagger} \qn \prefix{j}v}
\eeq

But

\beq
\prefix{j}w^{\dagger} ( E - \qn ) \prefix{j}v =
\left( w^{\dagger} - \prefix{j}y^{\dagger} \right) \left( E - \qn \right)
\left( v - \prefix{j}z \right) =
 \prefix{j}y^{\dagger}  \left( E - \qn \right) \prefix{j}z 
\eeq

\noindent because $w^{\dagger} (E  - \qn ) = 0$ and
$(E- \qn) v = 0$. Besides,

\beq
\prefix{j}w^{\dagger} \qn \prefix{j}v =
( \prefix{j}w^{\dagger} E v) (w^{\dagger}E  \prefix{j}v ) =
\left( \prefix{j}w^{\dagger} E \prefix{j}v \right) ^2
\eeq

So

\beq
\prefix{j}\rho^{-1} = \frac{\prefix{j}y^{\dagger}  \left( E - \qn \right) \prefix{j}z}
                                             {\left( \prefix{j}w^{\dagger} E \prefix{j}v \right) ^2}
= \prefix{j}\tilde{y}^{\dagger}\prefix{j}\tilde{z}
\label{roti}
\eeq

(Remenber that $(E- \qn)^2 = E - \qn$). Convergence results will
be proved from formulae (\ref{zti}), (\ref{yti}) and (\ref{roti}).

Some bounds will be derived. The 2-norm,
denoted as $\| \bullet \|$, will be
used in the sequel. Firstly, note that

\beq
\| E - \qn \| \leq 1
\eeq

In fact, let us consider an arbitrary vector
$\bar{v}$. As $E$ is a projection matrix

\beq
\| E \bar{v} \| \leq \| \bar{v} \|
\eeq

On the other hand, $\bar{v}$ can be decomposed
in a component lying on the eigenvector $v$ and
a perpendicular component $v_\bot$:

\beq
\bar{v} = a v + b v_\bot
\eeq

So

\beq
\| (E- \qn ) \bar{v} \| = \| a (E - \qn) v + b (E- \qn) v_{\bot} \| =
\| b E v_{\bot} \| \leq \| b v_{\bot} \| \leq \| \bar{v} \|
\eeq

Let us assume that

\beq
\| \prefix{j}\tilde{z} \| \leq \delta_z \;\;,\;\; \| \prefix{j}\tilde{y} \| \leq \delta_y
\eeq

Then 

\beq
\| \frac{(E - \qn)\prefix{j}z}{\prefix{j}w^{\dagger}E\prefix{j}v} \| \leq \delta_z \;\;,\;\; 
\| \frac{\prefix{j}y^\dagger (E - \qn)}{\prefix{j}w^{\dagger}E\prefix{j}v} \| \leq \delta_y
\eeq

As $\| E - \qn \| \leq 1$,

\begin{eqnarray}
\| \prefix{j} z \| = \| v - \prefix{j}v \| & \leq & \delta_z \| \prefix{j}w^{\dagger} E \prefix{j}v \| \label{zn1} \\
\| \prefix{j} y \| = \| w - \prefix{j}w \| & \leq & \delta_y \| \prefix{j}w^{\dagger} E \prefix{j}v \label{yn1} \| 
\end{eqnarray}

As, because of (\ref{conper}),

\beq 
1 = w^{\dagger}Ev = (\prefix{j}w^{\dagger} + \prefix{j}y^{\dagger}) E
                                      (\prefix{j}v + \prefix{j}z) = 
\prefix{j}w^{\dagger} E \prefix{j}v + \prefix{j}y^{\dagger} E \prefix{j}z
\eeq 

\noindent it is fullfiled

\beq
\| \prefix{j}w^{\dagger} E \prefix{j}v  \| \leq 1
+ \| \prefix{j}y^{\dagger} E \prefix{j}z \|
\eeq

Therefore, from (\ref{zn1}) and (\ref{yn1}),

\beq
\| \prefix{j}w^{\dagger} E \prefix{j}v  \| \leq 
1 + \delta_z \delta_y \| \prefix{j}w^{\dagger} E \prefix{j}v  \|^2
\eeq

\noindent so

\beq
\| \prefix{j}w^{\dagger} E \prefix{j}v  \| \leq 
\frac{1 + \sqrt{1 - 4 \delta_z \delta_y}}{2} \leq 1
\eeq

And (\ref{zn1}) and (\ref{yn1}) yield

\begin{eqnarray}
\| \prefix{j} z \| = \| v - \prefix{j}v \| & \leq & \delta_z \label{zn2} \\
\| \prefix{j} y \| = \| w - \prefix{j}w \| & \leq & \delta_y \label{yn2} 
\end{eqnarray}

So

\begin{eqnarray}
\| E \prefix{j}v \| & \leq & \| E v \| + \delta_z \\
\| \prefix{j}w^{\dagger} E \| & \leq & \|w^{\dagger}E \| + \delta_y
\end{eqnarray}

Up to now, eigenvectors $v$ and $w$ are arbitrary subject to the condition
$w^{\dagger}Ev=1$. In the sequel, they will be chosen such that

\beq
K = \| Ev \| = \| Ew^{\dagger}\|\;\;,\; w^{\dagger}Ev = 1
\eeq

Then,

\begin{eqnarray}
\qn - \prefix{j}\qn & = & Evw^{\dagger}E - E\prefix{j}v\prefix{j}w^{\dagger}E \nonumber \\
& = &
 E\left(v - \prefix{j}v \right) w^{\dagger}E +
E\prefix{j}v \left( w^{\dagger} - \prefix{j}w^{\dagger} \right)E \nonumber \\
& = &
E\prefix{j}z w^{\dagger}E + E \prefix{j}v \prefix{j}y E
\end{eqnarray}

So, being $m$ the dimension of $\qn$,

\begin{eqnarray}
\| \qn - \prefix{j}\qn \| & = &
\| E\prefix{j}z w^{\dagger}E + E \prefix{j}v \prefix{j}y E \| \nonumber \\
& \leq &
\| E\prefix{j}z w^{\dagger}E\| + \|E \prefix{j}v \prefix{j}y E \| \nonumber \\
& \leq &
m^2 \left[ 
\| E\prefix{j}z\| \|w^{\dagger}E\| + \|E \prefix{j}v\|  \|\prefix{j}y E \| 
\right] \nonumber \\
& \leq &
m^2 \left[ \delta_z K + (K + \delta_z) \delta_y \right] \nonumber \\
& = &
m^2 \left[  (\delta_y + \delta_z) K + \delta_z \delta_y \right] 
\label{bouq}
\end{eqnarray}

Let us consider the definition of the matrix $R_z(\prefix{j}\qn, \prefix{j}\epsilon )$.
It is clear, because of continuity with respect to $\prefix{j}\epsilon$, that

\beq
\exists \delta_{\epsilon 1} > 0 \backslash \| \prefix{j}\epsilon \| <  \delta_{\epsilon 1} \Rightarrow 
\| R_z(\prefix{j}\qn, \prefix{j}\epsilon) \| \leq \sqrt{2}
\| E \prefix{j}\pn N(\lambda ,\prefix{j}\qn) \prefix{j}\pn E \| \| \prefix{j}\epsilon \| 
\eeq

But

\begin{eqnarray}
\| E \prefix{j}\pn N(\lambda ,\prefix{j}\qn) \prefix{j}\pn E \| &
\leq &
\| E \prefix{j}\pn\| \| N(\lambda ,\prefix{j}\qn)\| \| \prefix{j}\pn E \| \nonumber \\
& \leq &
\| E - \prefix{j}\qn\| \| N(\lambda ,\prefix{j}\qn)\| \| E - \prefix{j}\qn  \| \nonumber \\
& \leq &
\| N(\lambda ,\prefix{j}\qn)\| 
\end{eqnarray}

So

\beq
\exists \delta_{\epsilon 1} > 0 \backslash \| \prefix{j}\epsilon \| <  \delta_{\epsilon 1} \Rightarrow 
\| R_z(\prefix{j}\qn, \prefix{j}\epsilon ) \| \leq \sqrt{2}
\| N(\lambda ,\prefix{j}\qn)  \| \| \prefix{j}\epsilon \| 
\eeq

Because of (\ref{bouq}), $\qn$ and $\prefix{j}\qn$ are close whenever
$\prefix{j}\tilde{z}$ and $\prefix{j}\tilde{y}$ are small. Therefore,

\beq
\exists \delta_{z1} > 0 \backslash \| \prefix{j}\tilde{z} \| <  \delta_{z1}  \wedge 
\| \prefix{j}\tilde{y} \| <  \delta_{z1} \Rightarrow 
\| N( \lambda, \prefix{j}\qn ) \| \leq \sqrt{2}   \| N( \lambda, \qn ) \|
\eeq

So,

\begin{eqnarray}
\exists \delta_{\epsilon 1} > 0 \wedge \exists \delta_{z1} > 0 
\backslash \| \prefix{j}\epsilon \| <  \delta_{\epsilon 1} \wedge
\| \prefix{j}\tilde{z} \| <  \delta_{z1}  \wedge 
\| \prefix{j}\tilde{y} \| <  \delta_{z1} & \Rightarrow  & \nonumber \\
\| R_z(\prefix{j}\qn, \prefix{j}\epsilon) \| \leq 2
\| N(\lambda ,\qn)  \| \| \prefix{j}\epsilon \| & &
\end{eqnarray}

Bounds on $N(\lambda ,\qn) $ will bw provided in the next
appendix. Besides, from (\ref{zti}) and  (\ref{yti}) 

\begin{eqnarray}
\| \prefix{j+1}\tilde{z} \| & \leq & \left[
               1 + 2 \| \prefix{j}\tilde{y}\dagger R_z(\prefix{j}\qn, \prefix{j}\epsilon) \prefix{j}z \| \right]
               \| E - \qn \| \| R_z(\prefix{j}\qn, \prefix{j}\epsilon) \| \| \prefix{j}\tilde{z} \| \\
\| \prefix{j+1}\tilde{y} \| & \leq & \left[
               1 + 2 \| \prefix{j}\tilde{y}\dagger R_z(\prefix{j}\qn, \prefix{j}\epsilon) \prefix{j}z \| \right]
               \| E - \qn \| \| R_z(\prefix{j}\qn, \prefix{j}\epsilon) \| \| \prefix{j}\tilde{y} \| 
\end{eqnarray}

\noindent whenever $\| \prefix{j}\tilde{y}\dagger R_z(\prefix{j}\qn, \prefix{j}\epsilon) \prefix{j}z \| 
< 1$. Therefore,

\begin{eqnarray}
\exists \delta_{\epsilon 1} > 0 \wedge \exists \delta_{z1} > 0 
\backslash \| \prefix{j}\epsilon \| <  \delta_{\epsilon 1} \wedge
\| \prefix{j}\tilde{z} \| <  \delta_{z1}  \wedge 
\| \prefix{j}\tilde{y} \| <  \delta_{z1} & \Rightarrow & \\
\left\{ \begin{array}{lcr}
\| \prefix{j+1}\tilde{z} \| & \leq & 2
 \left[
               1 + 4  \| N(\lambda ,\qn)  \| \delta_{\epsilon 1} \delta_{z1}^2
\right]
                \| N(\lambda ,\qn)  \| \| \prefix{j}\epsilon \|  \| \prefix{j}\tilde{z} \| \\
\| \prefix{j+1}\tilde{y} \| & \leq & 2
 \left[
               1 + 4  \| N(\lambda ,\qn)  \| \delta_{\epsilon 1} \delta_{z1}^2
\right]
                \| N(\lambda ,\qn)  \| \| \prefix{j}\epsilon \|  \| \prefix{j}\tilde{y} \| 
\end{array} \right.
& & \nonumber
\end{eqnarray}

Let us define

\beq
\delta_{z2} = \min \{ \delta_{z1}, \sqrt{\frac{1}{4 \| N(\lambda ,\qn)  \| \delta_{\epsilon 1}}}, \}
>  0
\eeq

Then,

\begin{eqnarray}
\exists \delta_{\epsilon 1} > 0 \wedge \exists \delta_{z2} > 0 
\backslash \| \prefix{j}\epsilon \| <  \delta_{\epsilon 1} \wedge
\| \prefix{j}\tilde{z} \| <  \delta_{z2}  \wedge 
\| \prefix{j}\tilde{y} \| <  \delta_{z2} & \Rightarrow & \label{primpar} \\
\left\{ \begin{array}{lcr}
\| \prefix{j+1}\tilde{z} \| & \leq &
     4          \| N(\lambda ,\qn)  \| \| \prefix{j}\epsilon \|  \| \prefix{j}\tilde{z} \| \\
\| \prefix{j+1}\tilde{y} \| & \leq &
4               \| N(\lambda ,\qn)  \| \| \prefix{j}\epsilon \|  \| \prefix{j}\tilde{y} \| 
\end{array} \right.
& & \nonumber
\end{eqnarray}

On the other hand, from (\ref{lineal}), it is clear that

\beq
\exists \delta_{\epsilon 2} > 0 \backslash \| \prefix{j} \epsilon \| <  \delta_{\epsilon 2}
\Rightarrow 
 \| \prefix{j+1} \epsilon \| \leq 2 \| \prefix{j}\rho^{-1}\|  \| \prefix{j} \epsilon \| 
 = 
 2 \| \prefix{j}\tilde{z} \|  \| \prefix{j}\tilde{y} \|   \| \prefix{j} \epsilon \| 
\eeq

\noindent because of (\ref{roti}). Let us define

\beq
\delta_{\epsilon 3} = \min \{ \delta_{\epsilon 1}, \delta_{\epsilon 2} \} > 0
\label{ter}
\eeq

Then, from (\ref{primpar}) and (\ref{ter}),

\begin{eqnarray}
\exists \delta_{\epsilon 3} > 0 \wedge \exists \delta_{z2} > 0 
\backslash \| \prefix{j}\epsilon \| <  \delta_{\epsilon 3} \wedge
\| \prefix{j}\tilde{z} \| <  \delta_{z2}  \wedge 
\| \prefix{j}\tilde{y} \| <  \delta_{z2} & \Rightarrow & \label{defini} \\
\left\{ \begin{array}{lcr}
\| \prefix{j+1}\tilde{z} \| & \leq &
     4          \| N(\lambda ,\qn)  \| \| \prefix{j}\epsilon \|  \| \prefix{j}\tilde{z} \| \\
\| \prefix{j+1}\tilde{y} \| & \leq &
4               \| N(\lambda ,\qn)  \| \| \prefix{j}\epsilon \|  \| \prefix{j}\tilde{y} \| \\
 \| \prefix{j+1} \epsilon \| & \leq &
2 \| \prefix{j}\tilde{z} \|  \| \prefix{j}\tilde{y} \|   \| \prefix{j} \epsilon \| 
\end{array} \right.
& & \nonumber
\end{eqnarray}

Let us now to define

\beq
\delta_{\epsilon 4} = \min \{ \delta_{\epsilon 3}, \frac{1}{8 \| N(\lambda ,\qn)  \| }  \} > 0
\;\;, \;\; 
\delta_{z3} = \min \{ \delta_{z2}, \frac{1}{2} \}
\eeq

Then

\begin{eqnarray}
\exists \delta_{\epsilon 4} > 0 \wedge \exists \delta_{z3} > 0 
\backslash \| \prefix{j}\epsilon \| <  \delta_{\epsilon 4} \wedge
\| \prefix{j}\tilde{z} \| <  \delta_{z3}  \wedge 
\| \prefix{j}\tilde{y} \| <  \delta_{z4} & \Rightarrow & \label{conver} \\
\left\{ \begin{array}{lcr}
\| \prefix{j+1}\tilde{z} \| & \leq &
   \frac{1}{2} \| \prefix{j}\tilde{z} \| \\
\| \prefix{j+1}\tilde{y} \| & \leq & \frac{1}{2} \| \prefix{j}\tilde{y} \| \\
 \| \prefix{j+1} \epsilon \| & \leq &
\frac{1}{2}| \prefix{j} \epsilon \| 
\end{array} \right.
& & \nonumber
\end{eqnarray}

\noindent which implies that the algorithm converges if
inizializated in a neighborhood of the solution
$(\lambda, v, w)$.  In fact, if $\| \prefix{j} z \| <
\| \prefix{j}w^{\dagger}E\prefix{j}v \| \delta_{z3} \|$
then, because of (\ref{E50})
$\| \prefix{j}\tilde{z}\| < \delta_{z3}$. There is an analogous
formula for $\| \prefix{j}\tilde{y}\|$. In that case,
\ref{conver} applies and the algorithm converges.

In order to asses the quotient convergence factor, let us
write \ref{defini} as 

\beq
\left[ \begin{array}{c}
\log  \| \prefix{j+1}\tilde{z} \| \\
\log  \| \prefix{j+1}\tilde{y} \| \\
\log  \| \prefix{j+1}\epsilon  \| 
\end{array} \right]
\leq
\left[ \begin{array}{ccc}
1 & 0 & 1 \\
0 & 1 & 1 \\
1 & 1 & 1
\end{array} \right]
\left[ \begin{array}{c}
\log  \| \prefix{j}\tilde{z} \| \\
\log  \| \prefix{j}\tilde{y} \| \\
\log  \| \prefix{j}\epsilon  \| 
\end{array} \right]
+
\left[ \begin{array}{c}
\log  4  \| N(\lambda ,\qn)  \|  \\
\log  4  \| N(\lambda ,\qn)  \|  \\
\log  2 
\end{array} \right]
\eeq

By substituting $\leq$ by $=$ a majorant sucession
is obtained. Its asymptotic behaviour is controlled
by the dominat eigenvalue of the state matrix, which
happens to be $1 + \sqrt{2}$. Therefore, it is expected
that, asymptotically

\beq
\log  \| \prefix{j+1}\epsilon  \| \leq (1+\sqrt{2}) \log  \| \prefix{j}\epsilon  \| 
\eeq

\noindent which completes the proof of the theorem.

\section{Some bounds on $N(\lambda, \qn$)}
\label{ape}

The purporse of this appendix is to provide some
bounds on the matrix $N(\lambda, \qn)$,
when $\lambda$ is an eigenvector of the pair $(E,A)$
and the matrix $\qn$ is formed from its right and left
eigenvectors: $\qn = Evw^{\dagger}E$.

 It shall be assumed that the pair $(E,A)$ is 
{\it solvable} \cite{Yip}, that is, the pencil
$\mu E - A$ is regular for all $\mu$ but
a finite number. Of course one of these
$\mu$ is the sought eigenvalue $\lambda$.
Moreover. it shall be assumed that
$\lambda$ is a single eigenvalue.

Let us assume that

\beq
N (\lambda, \qn) ^{-1}= \lambda E - A + A\qn
\eeq

\noindent is singular.  Then, there is
a vector $x$ such that

\beq
\left(  \lambda E - A + A\qn \right) x = 0
\eeq

However, it is possible to write

\beq
x = \alpha v + z \;\;, w^{\dagger}E z = 0\;\;, w^{\dagger}E  v = 1
\eeq

So $\qn z = 0, \qn v = Ev$. Then,

\beq
\alpha \lambda E v + \left(  \lambda E - A \right) z = 0
\label{pr1}
\eeq

Premultiplying by the left eigenvector $w$,

\beq
\alpha \lambda = 0
\eeq

Then, from (\ref{pr1}), it is obtained that
$ \left(  \lambda E - A \right) z = 0$. But this is
impossible, becuase it is assumed that
$\lambda$ is a single eigenvalue. Therefore,
the matrix $\lambda E - A + A\qn$ is regular
and $N (\lambda, \qn) $ is bounded.

A similar reasoning can be done
if 

\beq
N (\lambda, \qn) ^{-1}= \lambda E - A + \qn A
\eeq

Lastly, let us assume that

\beq
N (\lambda, \qn) ^{-1}= \lambda E - A + A\qn + \qn A
\eeq

Then, it is obtained that

\beq
\alpha \left( \lambda E  + \qn A \right) v +
 \left(  \lambda E - A  + \qn A \right) z = 0
\label{pr2}
\eeq

But $A v = \lambda Ev$, so

\beq
2 \alpha  \lambda E   v +
 \left(  \lambda E - A  + \qn A \right) z = 0
\label{pr3}
\eeq

Premultiplying by $w^{\dagger}E $, and taking into
account that $w^{\dagger} E \qn = w^{\dagger}E$, it is obtained that

\beq
 2 \alpha \lambda + \lambda w^{\dagger} E z =
 \alpha \lambda  = 0
\eeq

So, from (\ref{pr3}),

\beq
\left(  \lambda E - A  + \qn A \right) z = 0
\eeq

\noindent which implies that the
matrix $\lambda E - A  + \qn A$ is singular,
or that $z = 0$. As 
 it has been proved that $\lambda E - A  + \qn A$
is regular, it must be $z = 0$. Then, the singular
vector $x = v$. From
(\ref{pr3}), it must be $\lambda E v = 0$. 
So:

\begin{itemize}

\item If $\lambda = 0$, then $N(\lambda,\qn)$ is singular, and
          $v$ is a singular vector.

\item If $\lambda \neq 0$, then we must have $E v = 0$. But
          $Ev = \frac{1}{\lambda}Av \neq 0$, because if $Av = 0$, then
          $\lambda = 0$. So, the matrix $N(\lambda,\qn)$ is regular.

\end{itemize}

\section{Composite models}
\label{apef}

Let us define:

\beq
x = \left[ \begin{array}{c}
            x_M \\ x_I \\ x_O \\ x_A
           \end{array} \right]
\eeq

The system equations are

\begin{eqnarray}
\left[ \begin{array}{cccc}
\diag (E_1, E_2, \ldots, E_l) & 0 & 0 & 0 \\
0 & 0 & 0 & 0 \\
0 & 0 & 0 & 0 \\
0 & 0 & 0 & 0 
\end{array} \right]
\left[ \begin{array}{c}
            \dot{x}_M \\ \dot{x}_I \\ \dot{x}_O \\ \dot{x}_A
           \end{array} \right]
& = & \\
\left[ \begin{array}{cccc}
\diag (A_1, A_2, \ldots, A_l) & \diag (B_1, B_2, \ldots, B_l) & 0 & 0 \\
\diag (A_1, A_2, \ldots, A_l) & \diag (B_1, B_2, \ldots, B_l) &  I & 0 \\
0 & - J_{11} & I & -J_{12} \\
0 & - J_{21} & 0 & -J_{22} 
\end{array} \right]
\left[ \begin{array}{c}
            x_M \\ x_I \\ x_O \\ x_A
           \end{array} \right]
& & \nonumber
\end{eqnarray}

\noindent where $I$ denotes the identity matrix.

This linear dynamic system has 
eigenvectors which can be partitioned 
analogously to the variables and equations.
Specifically, it is fulfilled:

\begin{eqnarray}
\lambda E v_M & = & A v_M + B v_I \\
0 & = & C v_M + D v_I - v_O \\
0 & = & -J_{11} v_I + v_0 - J_{12}v_A \\
0 & = & -J_{21}v_I - J_{22}v_A \\
\lambda w_M^{\dagger}E & = & w_M^{\dagger} A + w_I^{\dagger}C \\
0 & = & w_M^{\dagger} B + w_I^{\dagger} D -
             \left( w_O^{\dagger} J_{11} + w_A^{\dagger} J_{21} \right) \\
0 & = & - w_I^{\dagger} + w_O^{\dagger} \\
0 & = & -w_O^{\dagger} J_{12} - w_A^{\dagger}J_{22}
\end{eqnarray}

By introducing the vector

\beq
\tilde{w}_O^{\dagger} = w_O^{\dagger}J_{11} + w_A^{\dagger}J_{21}
\eeq

\noindent and writing the previous equations in a subsystem basis, it is obtained:

\begin{eqnarray}
\lambda E_k v_{Mk} & = & A_k v_{Mk} + B_k v_{Ik}  \label{eigc1} \\
v_{Ok} & = & C_k v_{Mk} + D_k v_{Ik} \label{eigc2} \\
\left[ \begin{array}{cc}
J_{11} & J_{12} \\
J_{21} & J_{22}
\end{array} \right]
\left[ \begin{array}{c}
v_I \\ v_A
\end{array} \right] 
& = &
\left[ \begin{array}{c}
v_O \\ 0
\end{array} \right] 
\label{eigc3}  \\
\lambda w_{Mk}^{\dagger} E_k  & = & w_{Mk}^{\dagger} A_k + w_{Mk}^{\dagger} B_k   \label{eigc4} \\
\tilde{w}_{Ok}^{\dagger} & = & w_{Mk}^{\dagger} C_k  + w_{Ik}^{\dagger} D_k  \label{eigc5} \\
\left[ 
w_I^{\dagger} \; w_A^{\dagger}
 \right] 
\left[ \begin{array}{cc}
J_{11} & J_{12} \\
J_{21} & J_{22}
\end{array} \right]
& = &
\left[
\tilde{w}_O^{\dagger} \; 0
 \right]  
\label{eigc6} 
\end{eqnarray}

It is easy to check that

\beq
\cf^{\dagger} E \ce = I_m \Rightarrow \cf_{Mk}^{\dagger}E_k\ce_{Mk} = I_{mk}
\;\;\;\; \forall k
\eeq

\noindent where $mk$ is the number of columns of $\ce_{Mk}$.
Also, the formulae (\ref{deco1}-\ref{deco4}) can be written as:

\begin{eqnarray}
v_{Mk} & = & \ce_{Mk} \alpha_{Mk} + z_{Mk} \label{decoc1} \\
v_{Ik} & = & z_{Ik} \label{decoc2} \\
v_{Ok} & = & z_{Ok} \label{decoc3} \\
v_A & = & z_A \label{decoc3b} \\
\cf_{Mk}^{\dagger} E_k z_{Mk} & = & 0 \label{decoc4} \\
w_{Mk} & = & \cf_{Mk} \beta_{Mk} + y_{Mk} \label{decoc5} \\
w_{Ik} & = & y_{Ik} \label{decoc6} \\
\tilde{w}_{Ok} & = & y_{Ok} \label{decoc7} \\
w_A & = & y_A \label{decoc7b} \\
\ce_{Mk}^{\dagger} E_k y_{Mk} & = & 0 \label{decoc8} 
\end{eqnarray}

\noindent and

\beq
\alpha = \left[ \begin{array}{c}
            \alpha_{M1} \\ \alpha_{M2} \\ \vdots \\ \alpha_{Ml}
           \end{array} \right]
\;\;
\beta = \left[ \begin{array}{c}
            \beta_{M1} \\ \beta_{M2} \\ \vdots \\ \beta_{Ml}
           \end{array} \right]
\eeq

From (\ref{eigc1}) and (\ref{decoc1}-\ref{decoc3}):

\beq
\lambda E_k \ce_{Mk} \alpha_{Mk} + \lambda E_k z_{Mk} =
A_k \ce_{Mk} \alpha_{Mk} + A_k z_{Mk} + B_k z_{Ik}
\label{basc1}
\eeq

Premultiplying by $\cf_{Mk}^{\dagger}$, and taking into
account  (\ref{decoc4}):

\beq
\lambda \alpha_{Mk} = \cf^{\dagger}_{Mk}A_k\ce_{Mk} \alpha_{Mk} +
   \cf^{\dagger}_{Mk}A_k z_{Mk} + \cf^{\dagger}_{Mk} B_k z_{Ik}
\label{basc11}
\eeq

Let us define the projections $\pn_k$ y $\qn_k$ by:

\beq
\pn_k = I_{mk} - E_k\ce_{Mk}\cf^{\dagger}_{Mk}E_k =
               I_{mk} - \qn_k
\eeq

It is easy to check:

\begin{eqnarray}
\pn_k E_k \ce_{Mk} & = & 0 \\
\pn_k E_k z_{Mk} & = & E_k z_{Mk}
\end{eqnarray}

Taking the above equations into account, and
premultiplying (\ref{basc1}) by $\pn_k$:

\beq
\lambda E_k z_{Mk} = \pn_k A_k \ce_{Mk} \alpha_{Mk} +
                                         \pn_k A_k z_{Mk} +
                                        \pn_k B_k z_{Ik}
\eeq

So

\beq
\left( \lambda E_k - A_k + \qn_k A_k \right) z_{Mk} =
\pn_k A_k \ce_{Mk} \alpha_{Mk} + \pn_k B_k z_{Ik}
\eeq

Taking into account that $\pn_k z_{Mk} = z_{Mk}$
and $\qn_k z_{Mk} = 0$, the above equation yields:

\beq
z_{Mk} =
\pn_k 
\left( \lambda E_k - A_k +[ \qn_k, A_k]_+ \right)^{-1} \pn_k
\left\{ A_k \ce_{Mk} \alpha_{Mk} +
B_k z_{Ik} \right\}
\label{zMk}
\eeq

After substitution in equation 
(\ref{basc11}) it yields:

\beq
\lambda \alpha_{Mk} = \left(
A_{rk} + H_{Ak}(\lambda) \right)\alpha_{Mk}
          + \left( 
B_{rk} + H_{Bk}(\lambda) \right) z_{Ik}
\label{eqc1}
\eeq

\noindent where the following matrices have been
defined:

\begin{eqnarray}
A_{rk} & = & \cf^{\dagger}_{Mk}A_k\ce_{Mk} \\
B_{rk} & = & \cf^{\dagger}_{Mk} B_k \\
H_{Ak} (\lambda) & = &
\cf^{\dagger}_{Mk} A_k \pn_k
\left( \lambda E_k - A_k + [\qn_k, A_k]_+ \right)^{-1} \pn_k
A_k \ce_{Mk} \\
H_{Bk}(\lambda) & = &
\cf^{\dagger}_{Mk} A_k \pn_k
\left( \lambda E_k - A_k + [\qn_k, A_k]_+ \right)^{-1} \pn_k
B_k
\end{eqnarray}

On the other hand, from equation (\ref{compo2}),
it is obtained:

\beq
z_{Ok} = C_k \ce_{Mk} \alpha_{Mk} + C_k z_{Mk} + D_k z_{Ik}
\eeq

After substitution of (\ref{zMk}) it is obtained:

\beq
z_{Ok} = \left(C_{rk} + H_{Ck} \right) \alpha_{Mk} +
                 \left( D_k + H_{Dk} \right) z_{Ik}
\label{eqc2}
\eeq

\noindent where the following matrices have been defined:

\begin{eqnarray}
C_{rk} & = & C_k \ce_{Mk} \\
H_{Ck} (\lambda) & = &
C_k \pn_k
\left( \lambda E_k - A_k + [\qn_k, A_k]_+ \right)^{-1} \pn_k
A_k \ce_{Mk} \\
H_{Dk}(\lambda) & = &
C_k
\pn_k
\left( \lambda E_k - A_k + [\qn_k, A_k]_+ \right)^{-1} \pn_k
B_k
\end{eqnarray}

Let us define the matrices:

\begin{eqnarray}
A_r & = & {\rm diag} (A_{r1} \ldots A_{rl}) \\
B_r & = & {\rm diag} (B_{r1} \ldots B_{rl}) \\
C_r & = & {\rm diag} (C_{r1} \ldots C_{rl}) \\
D & = & {\rm diag} (D_{1} \ldots D_{l}) \\
H_{A} & = & {\rm diag} (H_{A1} \ldots H_{Al}) \\
H_{B} & = & {\rm diag} (H_{B1} \ldots H_{Bl}) \\
H_{C} & = & {\rm diag} (H_{C1} \ldots H_{Cl}) \\
H_{D} & = & {\rm diag} (H_{D1} \ldots H_{Dl}) 
\end{eqnarray}

It is easy to check that $A_r = \cf^{\dagger}A\ce$, as computed
according the general formula. 
Therefore, the proposed
notation is consistent.
Note also that all these matrices are diagonal-block matrices,
which eases its computation.
Then, equations 
(\ref{eqc1},\ref{eqc2}) can be written as:

\begin{eqnarray}
\lambda \alpha & = & \left( A_r + H_A(\lambda) \right) \alpha +
                                        \left( B_r + H_B(\lambda) \right) z_I \label{casf1} \\
z_O                       & = & \left( C_r + H_C(\lambda) \right) \alpha +
                                        \left( D     + H_D(\lambda) \right) z_I \label{casf2} 
\end{eqnarray}

By using  equations (\ref{eigc3}) and (\ref{decoc3b}), it is obtained:

\begin{eqnarray}
\lambda \alpha & = &
\left\{ 
A_r +
H_A(\lambda) + \right.   \nonumber \\
& &
\left.
\left[
\left( B_r + H_B(\lambda) \right) \; 0 \right]
\left[ \begin{array}{cc}
J_{11}-(D+H_D(\lambda)) & J_{12} \\
J_{21} & J_{22}
\end{array} \right]^{-1}
\left[ \begin{array}{c} 
C_r + H_C(\lambda) \\ 0 
\end{array} \right]
\right\}
\alpha 
\end{eqnarray}

Therefore

\beq
H(\lambda) =
H_A(\lambda) +
\left[
\left( B_r + H_B(\lambda) \right) \; 0 \right]
\left[ \begin{array}{cc}
J_{11}-(D+H_D(\lambda)) & J_{12} \\
J_{21} & J_{22}
\end{array} \right]^{-1}
\left[ \begin{array}{c} 
C_r + H_C(\lambda) \\ 0 
\end{array} \right]
\eeq

\section{Proof of theorem \ref{teof}}
\label{apeg}

Let us assume that

\beq
\left[ A - \prefix{j-1}\lambda \left(E - \qn \right) \right] \prefix{j}v =
\prefix{j}\lambda \qn \prefix{j}v
\eeq

But

\beq
\qn \prefix{j}v = (E \ce \cf^{\dagger}E) (\ce \prefix{j}\alpha + \prefix{j}z) =
E \ce \prefix{j}\alpha
\eeq

So

 \beq
\left[ A - \prefix{j-1}\lambda \left( E - \qn \right) \right] (\ce \prefix{j}\alpha + \prefix{j}z) =
\prefix{j}\lambda E \ce \prefix{j}\alpha \label{apeg1b}
\eeq

Premultiplying by $\cf^{\dagger}$, and taking into account 

\beq
\cf^{\dagger}  \left( E - \qn \right) =
\cf^{\dagger}E - \cf^{\dagger}E E \ce \cf^{\dagger}E = 0
\eeq

\noindent it yields

\beq
\cf^{\dagger} A \ce \prefix{j}\alpha + \cf^{\dagger} A \prefix{j}z =
\prefix{j}\lambda \prefix{j}\alpha \label{apeg10}
\eeq

On the other hand, as

\beq
  \left( E - \qn \right) \ce =
E\ce -E \ce \cf^{\dagger}EE \ce = 0
\eeq

\noindent equation (\ref{apeg1b}) yields

 \beq
A \ce \prefix{j}\alpha +
\left[ A - \prefix{j-1}\lambda \left( E - \qn \right) \right] \prefix{j}z =
\prefix{j}\lambda E \ce \prefix{j}\alpha
\eeq

Premultiplying by $\pn$, and taking into account that
$\pn (E-\qn) = E - \qn$, that $\pn E \ce = (E - \qn)\ce = 0$,
and that $\qn \prefix{j}z = 0$,

\beq
\pn A \ce \prefix{j}\alpha + \pn A \prefix{j}z =
\prefix{j-1}\lambda E \prefix{j}z 
\eeq

This is, esentially, formula (\ref{A7}). So, it is possible to conclude
(\ref{A10})

\beq
\prefix{j}z = \pn \left\{ \prefix{j-1}\lambda E - A + \left[A,\qn\right]_+ \right\}^{-1}
                       \pn A \ce \prefix{j}\alpha
\eeq

\noindent and, by substituting in (\ref{apeg10}),

\beq
\prefix{j}\lambda \prefix{j}\alpha = \left( A_{rr} + H(\prefix{j-1}\lambda) \right)
\prefix{j}\alpha 
\eeq

Analogous reasoning can be done by using
the left eigenvector, which proves the theorem.
\section{Proof of formula (\ref{fN})}
\label{apeh}

From the $\prefix{j}\overline{V}$ definition (\ref{defV}):

\beq
\prefix{j}\overline{V} =
\left[ A - \prefix{j-1} \lambda \left( E - E\ce \cf^{\dagger}E \right) \right]^{-1}
E \ce
\eeq

By using the Shermann-Morrison lemma:

\begin{eqnarray}
\prefix{j}\overline{V} & = &
\left[ \left( A - \prefix{j-1} \lambda  E \right) - 
\prefix{j-1} \lambda  E\ce \cf^{\dagger}E  \right]^{-1}
E \ce \nonumber \\
& = &
\left\{
 \left( A - \prefix{j-1} \lambda  E \right)^{-1} -
\prefix{j-1}\lambda  \left( A - \prefix{j-1} \lambda  E \right)^{-1}
\nonumber 
\right. \\ 
& & \left.
\left[ I_n + \cf^{\dagger}E
 \left( A - \prefix{j-1} \lambda  E \right)^{-1}E\ce \prefix{j-1}\lambda \right]^{-1}
\cf^{\dagger} E  \left( A - \prefix{j-1} \lambda  E \right)^{-1} \right\} E \ce \nonumber \\
& = &
\left\{ \prefix{j}V - \prefix{j-1}\lambda \prefix{j}V \left[ I_n +
\prefix{j-1}\lambda \cf^{\dagger}E \prefix{j}V \right]^{-1}
 \cf^{\dagger}E \prefix{j}V \right\}
\end{eqnarray}

Taking into account that $\prefix{j}{\cal N}^{-1} = \cf^{\dagger}E\prefix{j}\overline{V}$
and that
$\prefix{j}{\cal M} = \cf^{\dagger}E\prefix{j}V$,

\begin{eqnarray}
\prefix{j}{\cal N}^{-1} & = & \cf^{\dagger}E\prefix{j}\overline{V} \nonumber \\
& = &
\prefix{j}{\cal M} - \prefix{j-1}\lambda \prefix{j}{\cal M}  \left[
I_n + \prefix{j-1}\lambda \prefix{j}{\cal M} \right]^{-1} \prefix{j}{\cal M} \nonumber \\
& = &
\prefix{j}{\cal M}
\left[ I_n + \prefix{j-1}\lambda \prefix{j}{\cal M} \right]^{-1}
\left( 
\left[ I_n + \prefix{j-1}\lambda \prefix{j}{\cal M} \right] 
- \prefix{j-1}\lambda \prefix{j}{\cal M} \right) \nonumber \\
& = &
\prefix{j}{\cal M}
\left[ I_n + \prefix{j-1}\lambda \prefix{j}{\cal M} \right]^{-1}
\end{eqnarray}

So,

\beq
\prefix{j}{\cal N} =
\left[ I_n + \prefix{j-1}\lambda \prefix{j}{\cal M} \right] 
\prefix{j}{\cal M}^{-1} =
\prefix{j}{\cal M}^{-1} + \prefix{j-1}\lambda I_n
\eeq

\section{Proof of theorem \ref{teoHN}}
\label{apei}

Let us consider the equation

\beq
\left[ A - \lambda (E-\qn) \right]v = E\ce \alpha '
\label{apei1}
\eeq

\noindent where $\alpha '$ is an arbitrary vector. Without loss
of generality, it can be written:

\beq
v = \ce \alpha + z \;\;, \cf^{\dagger}Ez = 0
\eeq

Premultiplying ($\ref{apei1}$) by $\cf^{\dagger}$,

\beq
\cf^{\dagger} \left[ A - \lambda (E-\qn) \right]v = \alpha '
\eeq

As $\cf^{\dagger}(E- \qn) = 0$,

\beq
\cf^{\dagger}A \ce \alpha + \cf^{\dagger}A z = \alpha '
\eeq

On the other hand, from (\ref{apei1}), as $(E-\qn)\ce = 0$,

\beq
Av - \lambda (E-\qn) z = E \ce \alpha '
\label{apei2}
\eeq

So,

\beq
A \ce \alpha + \left[ A - \lambda(E-\qn) \right] z = E \ce \alpha '
\eeq

As $\pn (E - \qn ) = E - \qn$, $\pn E \ce = 0$ and $\qn z = 0$,
premultiplying by $\pn$,

\beq
\pn A \ce \alpha + \pn A z - \lambda E z = 0
\eeq

And, as $\qn z = 0$ and $\pn z = z$,

\beq
z = \pn \left[ \lambda E - A + [A,Q]_+ \right]^{-1} \pn A \ce \alpha 
\eeq

So, from (\ref{apei2}),

\begin{eqnarray}
\alpha ' & = &
\cf^{\dagger}A \ce \alpha +
\cf^{\dagger} A  \pn \left[ \lambda E - A + [A,Q]_+ \right]^{-1} \pn A \ce \alpha 
\nonumber \\
& = &
\left( A_{rr} + H(\lambda) \right) \alpha \label{apei3}
\end{eqnarray}

On the other hand, let us define

\beq
\overline{V} = \left[ A - \lambda (E-\qn) \right]^{-1}E\ce
\eeq

So,

\beq
v = \overline{V}\alpha ' = \ce \alpha + z
\eeq

Premultiplying by $\cf^{\dagger}E$, and taking into account 
$\cf^{\dagger}E\ce = I_n$,

\beq
\cf^{\dagger}E \overline{V} \alpha ' = \alpha
\eeq

So

\beq
\left(\cf^{\dagger}E \overline{V} \right)^{-1}\alpha  = {\cal N} \alpha = \alpha '
\label{apei4}
\eeq

Now, $\alpha '$ is arbitrary. So, as (\ref{apei3}) and (\ref{apei4})
hold for any $\alpha'$, the matrices must be equal:

\beq
{\cal N} = A_{rr} + H(\lambda)
\eeq

\noindent which proves the claim.


\begin{thebibliography}{99}

\bibitem{Perez}{\sc I. J. Perez-Arriaga, G. Verghese, F. Schweppe},
               {\em Selective Modal Analysis With Applications to Electric Power
                          Systems}, IEEE Trans. on PAS, Vol. PAS-101,
              No. 9, Sept. 1982.

\bibitem{Auto}{\sc I. J. Perez-Arriaga, G. C. Verghese, F. L. Pagola, F. C. Schweppe},
{\em Developments in Selective Modal Analysis of
Small-Signal Stability in Electric Power Systems}.
Automatica, Vol 26, No 2, pp 215-231, 1990.

\bibitem{TesisPerez}{\sc I. J. Perez-Arriaga},
                {\em Selective Modal Analysis With Applications to Electric Power
                          Systems}, Ph. D. Thesis, Electrical Engineering, M. I. T.,
                         June, 1981.

\bibitem{SMAx}{\sc R. Criado, J. Soto, J. Corera, L. Rouco, I. J. Perez-Arriaga},
                {\em SMAS3: A state-of-the-art computer package 
                  for analysis of small signal stability
                  in large electric power systems},
     CIGRE Study Committee 38: Colloquium on Power System Dynamic
     Performance. Florianopolis (Brasil), 22-23 Septiembre 1993.

\bibitem{Sancha}{\sc J. L. Sancha, I. J. Perez-Arriaga},
     {\em Selective Modal Analysis of Power System Oscillatory Instability},
     IEEE Transactions on Power Systems, Vol. PWRS-3, No. 2, May 1988, pp.
     429-438.        

\bibitem{Rouco}{\sc  L. Rouco, I. J. Perez-Arriaga},
                {\em Multi-area analysis of small signal stability in large electric power systems by SMA},
     IEEE Transactions on Power Systems, Vol. PWRS-8, No. 3, August 1993, pp.
     1257-1265.        

\bibitem{Yip}{\sc E. L. Yip, R. F. Sincovec},
               {\em Solvability, Controllabilty and Observability of Continous Descriptor
                         Systems}, IEEE Trans. on Automatic Control, Vol. AC-26, No. 3, June 1981.


\end{thebibliography}
\end{document}